\documentclass[a4paper]{amsart}

\usepackage{a4wide,color,eucal,enumerate}
\title[A Family of Nonlinear Fourth Order Equations]
{A Family of Nonlinear Fourth Order Equations \\ of Gradient Flow Type}
\date{\today}

\author{Daniel Matthes}
\address{Technische Universit\"at Wien. Wiedner Hauptstra\ss e 8 E 101 -- 1040 Wien, Austria.}
\email{matthes@asc.tuwien.ac.at }
\thanks{D.M. has been partially supported by the Deutsche Forschungsgemeinschaft, grant JU 359/7.
 He thanks the Department of Mathematics of the Universit\`a di Pavia,
 where part of this research has been carried out, for the kind hospitality.}
\author{Robert J. McCann}
\address{Department of Mathematics, University of Toronto.  Toronto Ontario M5S 2E4 Canada.}
\email{mccann@math.toronto.edu}
\thanks{R.J.M. has been partially supported by
US National Science Foundation grant DMS 0354729
and the Natural Sciences and Engineering
Research Council of Canada grants RGPIN 217006-03 and -08.}
\author{Giuseppe Savar\'e}
\address{Dipartimento di Matematica, Universit\`a di Pavia. Via Ferrata, 1 -- 27100 Pavia, Italy.}
\email{giuseppe.savare@unipv.it}
\thanks{G.S. has been partially supported by MIUR-PRIN'06 grant for the project
 ``Variational methods in optimal mass transportation and in geometric measure theory''.}

\newcommand{\FG}{G}
\newcommand{\restr}[1]{\lower3pt\hbox{$|_{#1}$}}
\newcommand{\Restr}[1]{\lower3pt\hbox{$\big|_{#1}$}}
\newcommand{\umin}{u_{\mathrm{min}}}
\newcommand{\eps}{\varepsilon}
\newcommand{\topref}[2]{\stackrel{\eqref{#1}}#2}
\newcommand{\rr}{{\mbox{\boldmath$r$}}}
\newcommand{\ttau}{{\mbox{\boldmath$\tau$}}}

\newcommand{\eeta}{{\mbox{\boldmath$\eta$}}}
\newcommand{\sttau}{{\mbox{\boldmath$\scriptstyle\tau$}}}
\newcommand{\seeta}{{\mbox{\boldmath$\scriptstyle\eta$}}}
\newcommand{\ssttau}{{\mbox{\boldmath$\scriptscriptstyle\tau$}}}
\newcommand{\ee}{{\mbox{\boldmath$e$}}}
\newcommand{\qq}{{\mbox{\boldmath$q$}}}
\newcommand{\yy}{{\mbox{\boldmath$y$}}}
\newcommand{\xxi}{{\mbox{\boldmath$\xi$}}}
\newcommand{\zzeta}{{\mbox{\boldmath$\zeta$}}}
\renewcommand{\d}{{\rm d}}
\newcommand{\dx}{\d x}
\newcommand{\Rd}{{\setR^d}}
\newcommand{\Nolin}[2]{{\nolin[#1;#2]}}
\newcommand{\Nolinl}[2]{{\nolinl[#1;#2]}}
\newcommand{\Mu}{M}
\newcommand{\sfc}{{\mathsf c}}
\newcommand{\sft}{{\mathsf t}}
\newcommand{\sfm}{{\mathfrak m}}
\newcommand{\sfC}{{\mathsf C}}
\newcommand{\sfa}{{\mathsf a}}
\newcommand{\sfb}{{\mathsf b}}
\newcommand{\sfA}{{\mathsf A}}
\newcommand{\Mom}[2]{\mathsf m_{#1}[#2]}
\newcommand{\QMom}[1]{\mathfrak m_2[#1]}
\newcommand{\Mass}[1]{\mathfrak m(#1)}
\newcommand{\scaling}{\delta}
\newcommand{\rmd}{\mathrm d}
\newcommand{\rmD}{\mathrm D}
\newcommand{\rmX}{{\mathrm X}}
\newcommand{\sgrp}{{\mathsf S}^\V}
\newcommand{\Sgrp}[2]{\sgrp_{#1}(#2)}
\newcommand{\Dil}[1]{\mathfrak d_{#1}}

\newcommand{\setR}{{\mathbf R}}
\newcommand{\setN}{{\mathbf N}}
\newcommand{\df}{\operatorname{D}}
\newcommand{\vcf}{\zzeta}
\newcommand{\flw}{{\mathbf X}}

\newcommand{\tmap}{{\mathbf t}}

\newcommand{\velo}{{\mathbf v}}

\newcommand{\dv}{\operatorname{div}}
\newcommand{\pp}{{\mathcal P}}
\newcommand{\nolin}{{\mathcal N}_\alpha}
\newcommand{\nolinl}{{\mathcal N}_{\alpha,\lambda}}
\newcommand{\id}{\operatorname{id}}
\newcommand{\w}{u}
\newcommand{\s}{\sigma}
\newcommand{\ww}{v}
\newcommand{\sw}{\rho}
\newcommand{\www}{w}
\newcommand{\rel}{{\mathcal H}}
\newcommand{\fish}{{F}}
\newcommand{\ent}{{H}}
\newcommand{\Ent}[2]{\ent_{#1}(#2)}
\newcommand{\Bent}[1]{H(#1)}
\newcommand{\Fish}[2]{\fish_{#1}(#2)}
\newcommand{\Entm}{{\mathcal H}}
\newcommand{\Fishm}{{\mathcal F}}

\newcommand{\lbg}[1]{{{\mathcal L}^{#1}}}
\newcommand{\lbgd}{{\lbg d}}
\newcommand{\blatt}{b}
\newcommand{\Blatt}{{\beta}}

\newcommand{\al}{{\alpha,\lambda}}

\newcommand{\dom}{{\rm Dom}}
\newcommand{\U}{{\mathcal U}}
\newcommand{\V}{{\mathcal V}}

\newtheorem{theorem}{Theorem}[section]
\newtheorem{remark}{Remark}[section]

\newtheorem{lemma}[theorem]{Lemma}
\newtheorem{corollary}[theorem]{Corollary}
\newtheorem{definition}[theorem]{Definition}

\begin{document}

\begin{abstract}
 Global existence and long-time 
 behavior
 of solutions to a family of nonlinear fourth order evolution equations on $\setR^d$ are studied.
 These equations constitute gradient flows for the perturbed information functionals
 \begin{align*}
   \Fish{\alpha,\lambda}\w = \frac1{2\alpha}
   \int_{\setR^d} \big|\df (\w^{\alpha})\big|^2 \,\dx + \frac\lambda2 \int_{\setR^d}|x|^2 \w\,\dx
 \end{align*}
 with respect to the $L^2$-Wasserstein metric.
 The value of $\alpha$ ranges from $\alpha=1/2$,
 corresponding to a simplified quantum drift diffusion model,
 to $\alpha=1$, corresponding to a thin film type equation.
\end{abstract}
\numberwithin{equation}{section}
\maketitle

\section{Introduction and Main Results}
This paper is concerned with non-negative solutions
to the following family of nonlinear fourth order parabolic problems in $(0,+\infty)\times \setR^d$,
\begin{align}
 \label{eq.main}
 & \partial_t\w = 
 - \dv\Big( \w \df\big[\w^{\alpha-1}\Delta \w^{\alpha}\big]\Big) + \lambda \dv( x \w )  = 0, \\
 \label{eq.ic}
 & \w(0,x) = \w_0(x) \geq 0, \quad \int_\Rd \big(1 + |x|^2 \big) \w_0(x)\, \dx <+\infty .
\end{align}
For the exponent $\alpha$ in \eqref{eq.main} we consider values $1/2\leq\alpha\leq1$,
and we assume $\lambda\geq0$.
%
This class of equations has recently been considered by
Denzler and McCann, who constructed special solutions \cite{DM}.
However, apart from the special cases $\alpha=1/2$ and $\alpha=1$,
which are described in some detail below,
analytical results for general solutions to \eqref{eq.main}
do not seem to have been available until now.
The current paper provides proofs
for the existence of weak solutions and their long-time behavior
under mild assumptions on $\w_0$.

Solutions to \eqref{eq.main} can be interpreted as nonlinear (confined) diffusion processes.
In fact, our results show remarkable similarities between \eqref{eq.main}
and the well-studied second-order (non-)linear Fokker Planck equations \cite{vazquez},
\begin{subequations}
 \label{eq.superdiff}
 \begin{align}
   \label{eq.slowdiff}
   \partial_t \ww = \Theta_\alpha\, \Delta
   (\ww^{\alpha+1/2})
   +  \Lambda_\al\, \dv(x\ww) \quad \mbox{if $\alpha>1/2$}, \\
   \label{eq:80}
   \partial_t \ww = \frac 12 \Delta
   \ww
   +  \Lambda_{1/2,\lambda} \dv(x\ww)
   \quad \mbox {if $\alpha=1/2$},
 \end{align}
\end{subequations}
with parameters defined by
\begin{align}
 \label{eq.slowdiffparams}
 \quad \Theta_{\alpha}:=\frac{\sqrt {2\alpha}}{2\alpha+1},
 \quad \Lambda_{\alpha,\lambda} := \sqrt{\lambda/\scaling_\alpha},
 \quad \scaling_\alpha:=(2\alpha-1)d+2.
\end{align}
%
The first similarity is that both equations \eqref{eq.main} and \eqref{eq.superdiff}
admit non-negative global weak solutions $\w$ and $\ww$, respectively.
The preservation of non-negativity by solutions to \eqref{eq.main} is a remarkable fact
in its own right,
since maximum principles do not generally apply to fourth order equations.
Moreover, the integral (total mass)
\begin{equation}
 \label{eq:81}
 \sfm=\Mass u:=\int_\Rd u\,\dx
\end{equation}
is preserved in time for both equations.
It is no loss of generality to assume unit mass $\sfm = 1$ in the following.

Further, as is characteristic for diffusion processes,
both solutions $\w$ and $\ww$ dissipate a variety of entropy functionals,
like the perturbed entropy $H_\al$ and perturbed information $F_\al$
introduced in \eqref{eq.info}--\eqref{eq.entropy1} below.
For $\lambda>0$, this leads to the convergence of $\w$ and $\ww$
to stationary solutions $\w_\infty$ and $\ww_\infty$ in the limit $t\to\infty$.
The attained profiles $\w_\infty$ and $\ww_\infty$
are independent of the initial conditions $\w_0$ and $\ww_0$, respectively.

In fact, the most striking similarity of \eqref{eq.main} and \eqref{eq.superdiff},
first observed in \cite{DM},
is that the stationary solutions $\w_\infty$ and $\ww_\infty$ are identical.
They are Barenblatt profiles or Gaussians, respectively,
\begin{subequations}
 \begin{align}
   \label{eq.barenblatt}
   \blatt_{\al;\sfm}(x) &=  \big( \sfa - \sfb |x|^2 \big)_+^{1/(\alpha-1/2)} ,
   \quad  \sfb = \frac{\alpha-1/2}{\sqrt{2\alpha}} \Lambda_{\alpha,\lambda} \quad & \mbox{if $\alpha>1/2$}, \\
   \label{eq.gauss}
   \blatt_{1/2,\lambda;\sfm}(x) &= \sfa\,\exp\big(-\Lambda_{1/2,\lambda}\,|x|^2
   \big)
   & \mbox{if $\alpha=1/2$}.
 \end{align}
\end{subequations}
Above, the positive parameter $\sfa$ is chosen to
adjust the mass $\sfm$ to unity.

One of the goals of this paper is to provide an explanation
of the aforementioned similarities between the fourth and the second order equations.
By exploiting the gradient flow structure of \eqref{eq.main},
these similarities turn out to reflect a natural correspondence
between entropy and information functionals.

\subsection*{Fourth order equations and gradient flows}

The family \eqref{eq.main} of nonlinear fourth order equations
possesses a gradient flow structure \cite{DM}.
This fact has been made use of by
Giacomelli and Otto \cite{ottox} \cite{GiOt} in the case $\alpha=1$,
and by Gianazza, Savar\'e and Toscani \cite{GST} in the case $\alpha=1/2$,
but it seems not to have been exploited previously for the intermediate range $1/2<\alpha<1$.

The suitable metric space is that of probability measures on $\setR^d$,
endowed with the $L^2$-Wasserstein distance.
Observe that equation \eqref{eq.main} is formally equivalent to
the canonical form of a Wasserstein gradient flow \cite[Ex. 11.1.2]{AGS}
\begin{equation}
 \label{eq:33}
 \partial_t\w + \dv\big(u \velo \big) = 0, \quad \velo = - \df \frac{\delta\fish_{\alpha,\lambda}(\w)}{\delta u},
\end{equation}
where $F_{\alpha,\lambda}$ is the perturbed 
{\em information functional} defined for smooth positive densities as
\begin{align}
 \label{eq.info}
 \Fish{\alpha,\lambda}\w :=& \frac1{2\alpha} \int_{\setR^d} \big|
 \df\big( \w^\alpha \big) \big|^2\,\dx +
 \frac\lambda2 \int_{\setR^d} |x|^2\w\,\dx
\end{align}
and $\delta \fish_\al/\delta u$ denotes its Eulerian first variation,
\begin{align*}
 \frac{\delta\fish_\al}{\delta u}=-\nabla\cdot(\alpha u^{2\alpha-2} \df u)+
 \alpha(\alpha-1)u^{2\alpha-3}\big|\df u\big|^2+\frac\lambda2|x|^2=
 -u^{\alpha-1}\Delta(u^\alpha)+\frac\lambda2|x|^2.
\end{align*}
We remark that the range $1/2\leq\alpha\leq1$ corresponds to the range
of convexity of $\fish_\al$ with respect to linear interpolation of measures.
On the other hand, we emphasize that even in this range of $\alpha$'s, 
the functionals $\fish_\al$ are not geodesically convex (= displacement convex);
see Carrillo and Slep\v{c}ev \cite{CS}.
The latter fact makes it impossible to apply the standard machinery
for metric gradient flows developed e.g.\ by Ambrosio, Gigli and Savar\'e \cite{AGS}
in a straight-forward manner.

The information establishes the link between \eqref{eq.main}
and the second order diffusion equation \eqref{eq.superdiff}.
Namely, consider in addition the perturbed {\em entropy functionals}
\begin{subequations}
 \begin{align}
   \label{eq.entropy0}
   \Ent{\alpha,\lambda}{\w} &= \frac{\Theta_\alpha}{\alpha-1/2}
   \Big( \int_{\setR^d} \w^{\alpha+1/2}\,\dx - \Mass {\w}^{\alpha+1/2} \Big)
   + \frac{\Lambda_{\alpha,\lambda}}2 \int_{\setR^d} |x|^2 \w \,\dx
   \quad\text{if $\alpha>1/2,$}\\
   \label{eq.entropy1}\Ent{1/2,\lambda}{\w} &= \frac12
   \Big( \int_{\setR^d} \w \log\w\,\dx - \Mass {\w}\log\big(\Mass\w\big)\Big)
   + \frac{\Lambda_{1/2,\lambda}}2 \int_{\setR^d} |x|^2 \w \,\dx \quad
   \text{if $\alpha=1/2$},
 \end{align}
\end{subequations}
which Otto showed generate \eqref{eq.superdiff}
as gradient flow with respect to the $L^2$-Wasserstein-metric \cite{otto}.
Define the associated {\em entropy production}
as the time derivative of $\ent_{\alpha,\lambda}$
along any sufficiently regular solution $\ww(t)$ to \eqref{eq.superdiff}.
The latter amounts to
\begin{align}
 \label{eq.specialform}
 - \frac{\d}{\d t}
 \Ent{\alpha,\lambda}{\ww(t)} =
 \Fish{\alpha,\lambda}{\ww(t)} -
 (\delta_\alpha-2)\Lambda\, \Ent{\alpha,\lambda}{\ww(t)}.
\end{align}
The special form of \eqref{eq.specialform},
as well as the choice of $\Lambda$ and $\delta$ in \eqref{eq.slowdiffparams},
is related to the rescaling properties of the unperturbed functionals $F_{\alpha,0}$ and $H_{\alpha,0}$
with respect to the ``mass preserving dilations''
\begin{equation}
 \label{eq:114}
 \Dil r u(\cdot):= r^{-d} u(r^{-1}\, \cdot)\quad
 r>0,
\end{equation}
which satisfy
\begin{equation}
 \label{eq:115}
 F_{\alpha,0}(\Dil r u)=r^{-\delta_\alpha} F(u),\qquad
 H_{\alpha,0}(\Dil r u)=r^{(-\delta_\alpha+2)/2} H_{\alpha,0}(u)\quad
 \forall\, r>0.
\end{equation}
The connection between \eqref{eq:115} and \eqref{eq.specialform} are clarified
in Remark \ref{rem:clarify}.

Relation \eqref{eq.specialform} is at the very basis of our investigations
of the long-time behavior of solutions $\w(t)$ to \eqref{eq.main}:
as $H_\al$ is a geodesically $\Lambda$-convex functional \cite{MC},
it generates a $\Lambda$-contractive gradient flow \eqref{eq.superdiff}, as in
Ambrosio, Gigli and Savar\'e \cite{AGS},
Carrillo McCann and Villani \cite{CarrilloMcCannVillani06}, and Sturm and
von Renesse \cite{SturmvonRenesse05} \cite{Sturm05}, after work of Otto \cite{otto}.
By means of \eqref{eq.specialform},
attraction towards the fixed point
is inherited by the gradient flow \eqref{eq.main},
regardless of the fact that its generating functional $F_\al$ is not geodesically convex.
The respective proof, in a very general setting,
is contained in section \ref{sec:basic estimate}.
%

\subsection*{The limiting cases: DLSS and thin film equation}

Another motivation for studying the family \eqref{eq.main}
originates from the equations obtained for $\alpha=1/2$ and $\alpha=1$.
When $\alpha=1/2$ and $\lambda=0$, the respective entropy \eqref{eq.entropy1}
turns into (half) the Boltzmann functional,
\begin{align}
 \label{eq.logentropy}
 H_{1/2,0}[\w] = \frac12 \Bent\w, \quad
 \Bent\w: = \int_{\setR^d} \w \log \Big(\frac{\w}{\Mass \w}\Big)\,\dx,
\end{align}
and the information coincides with the classical Fisher information,
\begin{align}
 \label{eq.fisher}
 \Fish{1/2,0}\w =
 \int_{\setR^d} \big|\df\sqrt{\w}\big|^2\,\dx =
  \frac14\int_\Rd \frac{|\df \w|^2}\w\,\dx.
\end{align}
The minimizers of the associated $\lambda$-perturbed functionals
are Gaussians \eqref{eq.gauss}.
For the respective unperturbed gradient flow,
one recovers the {\em Derrida-Lebowitz-Speer-Spohn equation}
\begin{align}
 \label{eq.qdbis}
 \partial_t \w + \dv\Big( \w \df \frac{\Delta\sqrt{\w}}{\sqrt{\w}} \Big) = 0.
\end{align}
Originally, those four authors derived \eqref{eq.qdbis}
in one spatial dimension on the half-line $\setR_+$
as a description for interface fluctuations in the Toom model \cite{DLSS1} \cite{DLSS2}.
Later, this equation was re-discovered in the mathematical theory of semiconductors.
In three spatial dimensions,
equation \eqref{eq.qdbis} constitutes the low-temperature, field-free limit
of a simplified quantum drift diffusion system for the electron density $\w$;
see \cite{JP} and references therein.
Existence of solutions to the initial value problem for \eqref{eq.qdbis}
has been studied by Bleher, Lebowitz and Speer \cite{BLS},
Gianazza, Savar\'e, and Toscani \cite{GST}, and
J\"ungel with Pinnau \cite{JP} and with Matthes \cite{JM}.
While still little is known about the qualitative properties of general solutions
(strict positivity is essentially an open problem, as is uniqueness),
a variety of results on the rates for relaxation to equilibrium
are available \cite{dlssreview}.

On the other hand,
the parameter $\alpha=1$ corresponds to 
\begin{align*}
 \Ent{1,0}\w = \frac{\sqrt{2}}{3} \Big( \int_{\setR^d} \w^{3/2}\,\dx - \Mass {\w}^{3/2} \Big), \quad
 \Fish{1,0}\w = \frac12 \int_{\setR^d} \big| \df \w \big|^2 \,\dx ,
\end{align*}
respectively, and their associated $\lambda$-perturbations.
Notice that $\Fish{1,0}\w$ is the Dirichlet energy of $\w$.
Both $\Ent{1,0}\cdot$ and $\Fish{1,0}\cdot$ are minimized by the Smyth-Hill profiles
\begin{displaymath}
 \, \blatt_{1,\lambda;\sfm} (x)
 =  \frac \lambda{8(d+2)}\big( \rho^2 - |x|^2
 \big)_+^2 ,\quad
 \text{$\rho>0$ being chosen to adjust the mass $\sfm$ to unity.}
\end{displaymath}
Equation \eqref{eq.main} turns into the following particular {\em thin film equation},
\begin{align}
 \label{eq.hs}
 \partial_t \w + \dv\big( \w \df\Delta \w) = 0.
\end{align}
The degenerate parabolic fourth order thin film (or lubrication) equations
of the more general form
\begin{align}
 \label{eq.tf}
 \partial_t \w + \dv\big( m(\w) \df\Delta \w) = 0
\end{align}
with an increasing, non-negative mobility function $m$
have been extensively studied in the mathematical literature
since the seminal paper of Bernis and Friedman \cite{BF}.
The theory of existence of weak solutions to \eqref{eq.tf} in dimensions $d=1,2,3$
with physically relevant mobility functions is quite complete,
see Bertsch, Dal Passo, Garcke and Gr\"un \cite{BDGG} and the references there.
Properties of the solution,
in particular the rate of convergence to the equilibrium state,
the spreading behavior of the support,
and waiting time phenomena
have been heavily investigated in a series of publications;
see e.g. \cite{DGG,G} and the review by Becker and Gr\"un \cite{BG}.
However, fundamental questions ---
like the rupture of film for specific mobility functions ---
still remain largely open, even in one spatial dimension.

The particular equation \eqref{eq.hs}, which is \eqref{eq.tf} with the linear mobility function $m(\xi)=\xi$,
plays a distinct role in the theory, at least in $d=1$ spatial dimension.
The equation generates the so-called {\em Hele-Shaw flow},
which describes the pinching of thin necks (with thickness $\w$) in a Hele-Shaw cell.
The mathematical literature on problems related to the Hele-Shaw flow is extensive;
an overview over classical results can be obtained e.g. in Myers \cite{hele}.
Analytical treatment in the spirit of our investigations is found in
Almgren, Bertozzi, and Brenner \cite{ABB}, Otto \cite{ottox},
and particularly in Carlen and  Ulusoy \cite{CU} and Carrillo and Toscani \cite{CT2}.

\subsection*{Main results.}
The results of this paper are twofold.
The first concerns the existence of global weak solutions
to the initial value problem \eqref{eq.main}\&\eqref{eq.ic}.

In order to formulate the existence theorem,
we introduce our notion of weak solution.
Observe that the non-linearity $\qq$
inside the first divergence operator of \eqref{eq.main}
can be formally rewritten as
\begin{equation}
 \label{eq:64} \qq
 =u\df\big(u^{\alpha-1}\Delta u^\alpha\big)=
 \df\big(u^\alpha\Delta u^\alpha)-
 u^{\alpha-1} (\Delta u^\alpha) \df u
 =\df\big(u^\alpha\Delta u^\alpha)-
 \alpha^{-1} (\Delta u^\alpha) \df u^\alpha.
\end{equation}
This expression has a (distributional)
meaning even if $\sigma=u^\alpha\in W^{2,2}_{\mathrm{loc}}(\Rd)$ only.
\begin{theorem}
 \label{thm.existence}
 Assume that the non-negative initial
 condition $\w_0\in L^1(\setR^d)$ satisfies
 \begin{equation}
   \label{eq:83}
   \int_\Rd |x|^2 \w_0(x)\,\d x<+\infty,\quad
   \Bent{ \w_0}=\int_\Rd \w_0(x)\log\big(\w_0(x)\big)\,\dx<+\infty.
 \end{equation}
 Then there exists a non-negative global solution
 $\w\in C^0([0,+\infty);L^{1}(\setR^d))$
 with $u(0)=u_0$ and $\w^\alpha\in L^2_{\mathrm{loc}}([0,+\infty);W^{2,2}(\Rd))$,
 satisfying the initial value problem \eqref{eq.main}\&\eqref{eq.ic}
 in the following sense
 \begin{align}
   \label{eq.weakform}
   \partial_t \w + \dv\Big(\df\big(u^\alpha\Delta u^\alpha)-
   \alpha^{-1} (\Delta u^\alpha) \df u^\alpha-\lambda u\,x\Big)=0,
 \end{align}
 i.e., for every test function
 $\zeta\in C^\infty_0((0,+\infty)\times\Rd)$
 one has
 \begin{equation}
   \label{eq:65}
   \int_0^{+\infty}\int_\Rd
   \big(-\partial_t\zeta+\lambda x\cdot\df \zeta\big)\w \,\dx\,\d t+
   \int_0^{+\infty}\int_\Rd \Delta u^\alpha \Big(u^\alpha\,
   \Delta \zeta +
   \alpha^{-1}\df u^\alpha\cdot\df\zeta\Big)\,\dx\,\d t=0.
 \end{equation}
\end{theorem}
\begin{remark}[Weaker formulation]
 \label{rem:weaker formulation}
 \upshape
 An even weaker form of \eqref{eq.main} exists,
 which just requires $\sigma=u^\alpha \in L^1_{\mathrm{loc}}((0,+\infty);W^{1,2}(\Rd))$.
 Indeed, observe that the vector $\qq$ defined in
 \eqref{eq:64} has components (repeated indices are summed)
 \begin{align*}
   \qq_i&=
   \partial_i\big(\sigma\partial^2_{jj}\sigma\big)-\alpha^{-1}
   (\partial_i\sigma)\partial^2_{jj}\sigma
   =
   \frac 12\partial_i\partial_{jj}^2(\sigma^2)
   -\partial_i|\partial_j\sigma|^2
   -\alpha^{-1}\partial_{j}\big(\partial_i\sigma\partial_j\sigma)
   +\alpha^{-1}(\partial^2_{ij}\sigma)\partial_j\sigma
   \\&=
   \frac 12\partial_i\partial_{jj}^2\big(\sigma^2)
   -(1-\frac 1{2\alpha})\partial_i|\partial_j\sigma|^2
   -\frac 1\alpha\partial_{j}\big(\partial_i\sigma\partial_j\sigma).
 \end{align*}
 The respective weaker formulation of \eqref{eq.weakform} reads for $\lambda=0$ as
 \begin{equation}
   \label{eq:66}
   \partial_t u+\frac 12\Delta^2(\sigma^2)-(1-\frac 1{2\alpha})\
   \Delta|\df \sigma|^2-\frac 1\alpha\sum_{ij}\partial^2_{ij}
   \big(\partial_i\sigma\partial_j\sigma)=0,\quad
   \sigma=u^\alpha.
 \end{equation}
 This corresponds to the integral condition
 \begin{equation}
   \label{eq:66bis}
    \int_0^{+\infty}\int_\Rd
   \big(-\partial_t\zeta+\lambda x\cdot\df \zeta\big)\w \,\dx\,\d t+
   \int_0^{+\infty} N_\alpha(u;\df \zeta) \, \d t=0
   \quad \forall\,\zeta\in C^\infty_0((0,+\infty);\Rd),
 \end{equation}
 upon defining for $u\in W^{1,1}(\Rd)$ with $\sigma=u^\alpha\in W^{1,2}(\Rd)$
 and $\zzeta=\df \zeta\in C^\infty_0(\Rd;\Rd)$
 \begin{equation}
   \label{eq:67}
   N_\alpha(u;\zzeta):=-\frac 1{2\alpha}\int_\Rd\Big(\alpha\,\df\sigma^2\cdot\df \dv\zzeta
   +(2\alpha-1)|\df\sigma|^2 \dv\zzeta+
   2 \df \zzeta \df \sigma\cdot\df\sigma\Big)\,\dx.
 \end{equation}
 It is not difficult to check that \eqref{eq:67} and \eqref{eq:65}
 are equivalent if $\sigma=u^\alpha\in L^1_{\mathrm{loc}}((0,+\infty);W^{2,2}(\Rd))$.
 Although the link of \eqref{eq.main} to \eqref{eq:66}
 is less evident than to \eqref{eq.weakform},
 the weaker form \eqref{eq:66} is well-suited
 to our analysis for several reasons:
 first, \eqref{eq:66} is the natural variational formulation
 of the gradient flow equation
 that is obtained by the method of Jordan, Kinderlehrer and Otto \cite{JKO},
 see Lemma \ref{lem.subdiff};
 second, the pointwise condition $\sigma=\w^\alpha\in W^{1,2}(\Rd)$ is satisfied
 as soon as $\fish_\al[\w]$ is finite;
 third, strong convergence of $\sigma=\w^\alpha$ in $L^1_{\mathrm{loc}}((0,+\infty);W^{1,2}(\Rd))$
 is sufficient to pass to the limit in the nonlinear operator $N_\alpha$.
 In particular, the last point is important for the passage
 to the continuous time limit in the minimizing movement scheme
 introduced in the next paragraph.
\end{remark}
The second result concerns the convergence of $\w$ in the long-time limit.
\begin{theorem}
 \label{thm.asymptotics}
 If $\lambda>0$ then
 the weak solutions $\w$ to the initial value problem \eqref{eq.main}\&\eqref{eq.ic}
 constructed as the small time-step limit of the minimizing movement variational scheme
 \eqref{eq:59intro} (which in particular satisfy the properties
 stated in Theorem \ref{thm.existence})
 converge exponentially fast towards the respective Barenblatt
 profile $\blatt_{\al;\sfm}$
 in $L^1(\setR^d)$,
 \begin{align}
   \label{eq.l1decay}
   \big\|\nobreakspace\w(t) - \blatt_{\al;\sfm} \big\|_{L^1(\Rd)}
   &\leq C\,\big( \ent_{\alpha,\lambda}(\w_0) - \ent_{\alpha,\lambda}(\blatt_{\al;\sfm}) \big)^{1/2}
   e^{-\lambda t} ,
 \end{align}
 where the constant $C>0$ only depends on
 $\alpha,\lambda,d,\sfm$.
 Moreover, the entropy and the information functionals
 $\ent_\al,\fish_{\alpha,\lambda}$
 decay along this solution at an exponential rate,
 \begin{align}
   \label{eq.energydecay}
   \ent_{\alpha,\lambda}(\w(t)) -
   \ent_{\alpha,\lambda}(\blatt_{\al;\sfm})
   &\leq \big( \ent_{\alpha,\lambda}(\w_0) - \ent_{\alpha,\lambda}(\blatt_{\al;\sfm}) \big)  e^{-2\lambda t},\\
   \fish_{\alpha,\lambda}(\w(t)) -
   \fish_{\alpha,\lambda}(\blatt_{\al;\sfm})
   &\leq \big( \fish_{\alpha,\lambda}(\w_0) - \fish_{\alpha,\lambda}(\blatt_{\al;\sfm}) \big)e^{-2\lambda t}.
   \label{eq.informationdecay}
 \end{align}
\end{theorem}
We refer to section \ref{sct.asymptotics} for more details,
and also an improved convergence result in $W^{1,2}(\Rd)$
for the ``thin film'' case $\alpha=1$.

In the unconfined case $\lambda=0$,
there exists no non-trivial stationary solution, because
$\w(t)$ tends to zero in $L^1_{loc}(\setR^d)$ as $t\to\infty$.
However, we are able to describe the {\em intermediate asymptotics}:
$\w$ approaches the self-similar spreading Barenblatt profile
\begin{equation}
 \label{eq:116}
 \blatt_{\alpha,0;\sfm}(t,\cdot):=R(t)^{-d} \blatt_{\alpha,1;\sfm}(R(t)^{-1}\,\cdot)=
 \Dil{R(t)}{\blatt_{\alpha,1;\sfm}},\quad
 R(t):=\big(1+(\delta+2)t\big)^{1/(\delta+2)},
\end{equation}
in $L^1(\setR^d)$ at an algebraic rate;
see Corollary \ref{cor.intermediate}.
Thanks to the scaling invariance \eqref{eq:115} of equation \eqref{eq.main},
the latter result can be recovered from \eqref{eq.l1decay}.
The corresponding rescaling argument is by now classical, see
e.g.\ Dolbeault and del Pino \cite{DdP}.
In Section \ref{sct.intermediate},
we will present a variational derivation of it,
that is entirely based on the invariance property \eqref{eq:115}.
In the particular case $\alpha=1$ of the thin film equation,
it leads the following result:
\begin{corollary}
 \label{cor:thinblatt}
 Assume $\alpha=1$ and $\lambda=0$.
 Then any weak solution $\w$
 to the initial value problem \eqref{eq.main}\&\eqref{eq.ic}
 constructed as limit of the minimizing movement variational scheme
 satisfies
 \begin{equation}
   \label{eq:117}
   \|\w(t,\cdot)-b_{1,0;\sfm}(t,\cdot)\|_{L^1(\Rd)}\leq C R(t)^{-1},\quad
   \|\df \w(t,\cdot)-\df b_{1,0;\sfm}(t,\cdot)\|_{L^2(\Rd)}\leq
   C\, R(t)^{-2-d/2}.
 \end{equation}
\end{corollary}
Variants of the estimates in Theorem \ref{thm.asymptotics} are already known is two special cases.
For the situation $\alpha=1/2$ of the Derrida-Lebowitz-Speer-Spohn equation \eqref{eq.qdbis},
results equivalent to those of Theorem \ref{thm.asymptotics}
have been recently obtained by Gianazza, Savar\'e and Toscani \cite{GST}.
Moreover, alternative estimates on entropy and information decay are available
for the Hele-Shaw flow \eqref{eq.hs} in dimension $d=1$,
\begin{align}
 \label{eq.hs1}
 \partial_t \w = - \big( \w \w_{xxx} \big)_x + (x\w)_x .
\end{align}
In \cite{CT2}, entropy decay for solutions to \eqref{eq.hs1}
--- at the same exponential rate ---
has been obtained for strong solutions $\w$.
Instead of resorting to the general formalism of Wasserstein gradient flows,
Carrillo and Toscani rewrite \eqref{eq.hs1} in an ingenious way,
so that all estimates can be made completely explicit.
As we show in Section \ref{sct.ct}, their appealing formal idea can in fact be extended
to the whole family of equations \eqref{eq.main} with $\alpha>1/2$,
and to arbitrary dimensions $d\geq1$.

Finally, we mention that the behavior of the perturbed information
\begin{align*}
 \fish_{1,1}[\w] = \frac 12 \int_\setR \big| \w_x(x) \big|^2\,\dx + \frac12 \int_\setR x^2\w(x)\,\dx
\end{align*}
along solutions to \eqref{eq.hs1} was studied recently.
By proving asymptotic equipartition of kinetic and potential energy,
an algebraic decay estimate was derived by Carlen and Ulusoy \cite{CU},
and subsequently improved by Carlen \cite{C},
leading to the same exponential equilibration behavior that we find.

\subsection*{The discrete variational scheme}

In order to prove Theorem \ref{thm.existence},
we closely follow the strategy introduced by
Jordan, Kinderlehrer and Otto \cite{JKO} in the framework of the Wasserstein space,
and further developed in the book of Ambrosio, Gigli and Savar\'e \cite{AGS}.
After mass renormalization $\sfm=1$,
we associate to the time dependent density function $\w(t)$
satisfying \eqref{eq.main} the (absolutely continuous) probability measures
$\mu(t):=\w(t,\cdot)\lbgd$ in $\Rd$.
The quadratic moment of $\mu(t)$ is finite for all times $t>0$
if it is initially at $t=0$.
Hence $\mu$ attains values in the Wasserstein space $\pp_2(\Rd)$,
metrized by the $L^2$-Wasserstein distance $W_2(\cdot,\cdot)$;
we recall the relevant definitions in section \ref{sct.preliminaries}.

A time-discrete approximation of solutions to \eqref{eq.main}\&\eqref{eq.ic}
is constructed by application of {De Giorgi's} ``minimizing movement'' scheme
to the information functional $\Fishm_\al[\mu]=\fish_\al(u)$
(defined on Borel probability measures).
We review this scheme below;
see section \ref{subsec:GF} for a more detailed description.

Let a partition $\mathcal P_\sttau$ of the time interval $[0,+\infty)$ be given,
\begin{equation}
 \label{eq:106}
 \mathcal P_\sttau:=\{0=t^0_\sttau<t^1_\sttau<\cdots<t^n_\sttau<\cdots\},\quad
 \tau_n:=t^n_\sttau-t^{n-1}_\sttau,\quad
 \lim_{n\to+\infty}t^n_\sttau=\sum_n\tau_n=+\infty.
\end{equation}
Associated to the step sizes $\ttau=(\tau_n)_{n\in\setN}$
and the initial measure $\Mu^0_\sttau=u_0\lbgd\in\pp_2(\setR^d)$,
we consider the sequence $(\Mu^n_\sttau)_{n\in \setN}$
recursively defined by solving
the following variational problem in $\pp_2(\setR^d)$:
\begin{equation}
 \label{eq:59intro}
 \text{find}\
 \Mu^n_\sttau\in \pp_2(\Rd)\
 \text{which minimizes the functional}\quad
 \Mu\mapsto \frac1{2\tau_n}W^2_2(\Mu^{n-1}_\sttau,\Mu) +
 \Fishm_\al[\Mu] .
\end{equation}
The measure $M^n_\sttau\in\pp_2$ at the $n$th time step ---
which serves as approximation of $\mu(t^n_\sttau)=u(t^n_\sttau)\lbgd$ ---
is thus determined as the minimizer of
a variational problem involving the Wasserstein distance and $\fish_\al$.

Following this scheme,
one obtains a family of piecewise constant approximating solutions
$M_\tau(t)\equiv M^n_\tau$ in each interval $(t^{n-1}_\sttau,t^n_\sttau]$,
satisfying a suitable discrete version of \eqref{eq:66}.
We prove that in the limit as $\sup_n\tau_n\downarrow0$,
one recovers a weak solution $\mu=u\lbgd$ of \eqref{eq.weakform}.
Since the functional $\fish_\al$ is not displacement convex,
we cannot follows the standard procedures developed in \cite{AGS}.
Strong a priori estimates are needed to pass to the limit
and to characterize the equation satisfied by the limit function.

\subsection*{Discrete dissipation estimates via auxiliary gradient flows}
The standard strategy for deriving these kinds of estimates is to get as
much information as possible from the discrete variational problems
\eqref{eq:59intro}, usually by studying the first variation of the
functional to be minimized
under suitable perturbations of the minimizer. The careful choice of such
perturbations therefore plays a crucial role.

In the present paper we propose a general strategy to find a
sufficiently wide family of perturbations and to obtain corresponding
discrete estimates, even in a general metric setting.
The underlying philosophy is simple:
we perturb the minimizer $M^n_\sttau$ of \eqref{eq:59intro}
by moving it along the gradient flows $\sgrp$ generated by other,
simpler functionals $\V$ defined on $\pp_2(\Rd)$.
In order to obtain useful information,
$\V$ should be $\kappa$-displacement convex for some $\kappa\in \setR$,
so that $\sgrp$ is $\kappa$-contracting
and can be characterized by a suitable variational inequality,
see Theorem \ref{thm.metricheat} and \eqref{eq:7}.
Assume $\kappa\geq0$ temporarily, for ease of presentation.
Assume further that the dissipation of $\Fishm_\al$
along the flow $\sgrp$,
\begin{equation}
 \label{eq:92}
 \df^\V\Fishm_\al[\mu]:=\limsup_{h\downarrow0}\frac{\Fishm_\al[\mu]-\Fishm_\al[\sgrp_h\mu]}h,
\end{equation}
is non-negative, possibly up to lower order terms.
Then, the discrete solution $M^n_\tau$ satisfies the a priori bound
\begin{equation}
 \label{eq:94}
 \V[M^n_\sttau]+\tau_n \df^\V\Fishm_\al[M^n_\sttau]\leq
 \V[M^{n-1}_\sttau]\quad
 \text{or, equivalently,}\quad
 \frac{\V[M^{n-1}_\sttau]-\V[M^n_\sttau]}{\tau_n}\geq \df^\V\Fishm_\al.
\end{equation}
See Theorem \ref{thm:metric_estimate} for the complete statement,
with $\kappa\in\setR$.
We shall refer to inequality \eqref{eq:94} ---
or rather to its generalized form \eqref{eq:13} ---
as the ``flow interchange estimate'' for the following reason.
The rate of dissipation $\tau^{-1}\big(\V[M^{n-1}_\tau]-\V[M^n_\tau]\big)$
of the functional $\V$ along the ``discrete flow''
generated by $\Fishm_\al$ through the \emph{minimizing movement scheme} \eqref{eq:59intro}
is bounded from below by the rate of dissipation $\df^\V\Fishm_\al$ of $\Fishm_\al$
along the continuous flow generated by $\V$.
It has the advantage of converting the differential estimate of $\Fishm_\al$
along a known evolution \eqref{eq:92} into a discrete estimate for the
approximation scheme of the flow for \eqref{eq.main}.

In the present paper we will repeatedly apply this principle
with various choices of $\V$:
\begin{enumerate}[\bf 1.]
\item
 {\em In Section \ref{subsec:first variation} and Lemma \ref{lem.velo}:}
 $\V[\mu]=\int_\Rd \zeta\,\d\mu$, for $\zeta\in C^\infty_0(\Rd)$,
 generating the flow associated to the linear transport equation
 \begin{equation}
   \label{eq:95}
   \partial_t v-\dv(v\df\zeta)=0.
 \end{equation}
 This first step (which goes back to the original ideas of
 \cite{JKO})
 shows that $M^n_\tau$ is a weak solution of a discrete version of
 \eqref{eq:66}.
\item {\em In Lemma \ref{le:noncompact} and equation \eqref{eq:84}:}
 $\V[\mu]=\QMom\mu :=\int_\Rd |x|^2\,\d\mu$, corresponding to the
 dilation equation
 \begin{equation}
   \label{eq:96}
   \partial_t v-\dv(v x)=0,\quad \sgrp_t v(x)=e^{dt}v(e^t x)=
   \Dil{e^{-t}}v\,(x).
 \end{equation}
 It provides a control on the quadratic moments of the solution in
 terms of the behavior of $\Fishm_\al$ with respect to the
 dilations group $\sgrp_t$, see \eqref{eq:115}.
\item {\em In Theorem \ref{thm:crucialapriori}:}
 $\V[\mu]=\rel[\mu]$, the logarithmic entropy functional from \eqref{eq.logentropy},
 generating the heat flow
 \begin{equation}
   \label{eq:97}
   \partial_tv-\Delta v=0.
 \end{equation}
 It provides the discrete prototype of the following a priori estimate
 for solutions $\w$ to \eqref{eq.main},
 \begin{equation}
   \label{eq:82}
   H(\w(T,\cdot))+\sfc_0 \int_0^T\int_\Rd |\df^2 \sigma(t,x)|^2\,\dx\d t\leq
   H(\w_0)\quad
   \forall\, T>0;\quad
   \sfc_0:=\frac {4d-1}{\alpha d(d+2)}.
 \end{equation}
 This estimate is crucial for the passage to the limit in the
 nonlinear discrete equation. With a different technique, an
 analogous idea has been exploited in \cite{GST}.
 In order to capture the full power of the estimate hidden in the
 term $\df^\V\Fishm_\al$, which controls the second derivatives
 of $u^\alpha$, we will use in a very careful way the integration
 by parts rule proved by Gianazza, Savar\'e and Toscani \cite{GST}
 and by J\"ungel and Matthes \cite{JM} (see also \cite{JM0} for
 the development of an algebraic technique
 for dealing with such formulas).
\item {\em In Theorem \ref{prp.prequilibrate}:}
 $\V[\mu]=\rel_\al[\mu]$, the power entropy \eqref{eq.entropy0},
 generating
 the porous medium flow \eqref{eq.slowdiff}.
 It provides all the asymptotic estimates on the exponential decay
 of $\rel_\al$ and $\Fishm_\al$ of Theorem \ref{thm.asymptotics}.
 This result is actually more subtle,
 since it also uses the particular relation \eqref{eq.specialform}
 between information function $\fish_{\alpha,\lambda}$ and the associated entropy $\ent_\al$.
\end{enumerate}

\subsection*{Open problems}
The restriction to the range $1/2\leq\alpha\leq1$ --- although natural for several reasons ---
is seemingly non-optimal with respect to 
where existence and equilibration results can be proven.
The authors suspect that all of the aforementioned results 
extend mutatis mutandis at least to
\begin{align}
  \label{eq.unnatural}
  \frac{(d-1)^2}{2d^2+1} < \alpha < \frac32 .
\end{align}
The reason is that on this wider range 
an a priori estimate of the type \eqref{eq:crucialapriori} continues to hold
for sufficiently regular solutions.

Still, even the wider range \eqref{eq.unnatural} constitutes a rather unexpected restriction.
In view of the properties of the intimately connected slow and fast diffusion flows,
the natural lower limit for $\alpha$ should be $\alpha=1/2-1/d$
(which is when the entropy $\ent_\al$ loses geodesical convexity),
while there should not exist any upper limit at all.

\subsection*{Plan of the paper.}
In Section \ref{sct.preliminaries} below
we collect some essential facts on the $L^2$-Wasserstein-metric,
its link with the entropy functionals (\ref{eq.entropy0},b) and their
gradient flows,
and the related properties of the information functional
\eqref{eq.info}.

As explained in the previous paragraph,
Section \ref{sec:basic estimate} is devoted to the core ideas of obtaining
dissipation estimates by using auxiliary gradient flows; its first
part presents the basic heuristics, deriving the
estimates in the elementary case of smooth flows in
the Euclidean space $\setR^m$. The general abstract results in metric
spaces
are developed in \ref{subsec:GF}.

The proof of the existence result Theorem \ref{thm.existence}, based
on the convergence of the minimizing movement scheme \eqref{eq:59intro},
is carried out in Section \ref{sct.existence}.
Besides the general ideas we have just tried to explain, the
crucial estimate and many technical steps are contained in Lemma
\ref{le:technical}.

The large-time asymptotic estimates from Theorem \ref{thm.asymptotics}
are derived in Section \ref{sct.asymptotics}, as a direct application
of the results of Section \ref{sec:basic estimate}.
Section \ref{sct.ct} and \ref{sct.intermediate} might be of independent interest.
In the first, a formal argument for estimate \eqref{eq.l1decay} is presented,
inspired by Carrillo and Toscani's approach to the special case $\alpha=d=1$ \cite{CT2}.
In the second, the intermediate asymptotics of the unconfined gradient flow,
i.e. equation \eqref{eq.main} with $\lambda=0$, are investigated.
The crucial tool is a rescaling property for the discrete solutions
of the minimizing movement scheme \eqref{eq:59intro}
for functionals satisfying \eqref{eq:115} under the action of the dilation
group $\mathfrak d_r$.

\section{Preliminaries}
\label{sct.preliminaries}

The framework of the proofs will be that of measures on $\setR^d$
and the $L^2$-Wasserstein distance.
For simplicity of notation, we restrict to probability measures from now on.
Due to the scaling invariance of \eqref{eq.main},
it costs little generality to consider only probability densities as solutions,
i.e. solutions of unit mass.
In fact, if $\w(t,x)$ is any solution to \eqref{eq.main} of mass
$\sfm>0$ in dimension $d\geq1$,
then the rescaled function
\begin{align*}
 \tilde u(t,\cdot) :=
 \sfm^{-1}\, Q^d\,  
 u\big( t, Q 
 \,\cdot)=\sfm^{-1} \Dil{Q^{-1}} u(t,\cdot), \quad \mbox{with}\quad
 Q:=\mathsf m^{(2\alpha-1)/4+(2\alpha-1)d},
\end{align*}
is another  solution to \eqref{eq.main}, and is of unit mass.

\subsection{Wasserstein distance}
\label{sct.wasserstein}
Denote by $\pp(\Rd)$ (resp.~$\pp_2(\Rd)$) the space of the Borel
probability measures on $\Rd$
(resp.\ with finite second moment $\QMom \mu=\int_\Rd |x|^2\,\d\mu$);
further, $\pp^r(\Rd)$ (resp.\ $\pp_2^r(\Rd)$) is the
subset of all the measures which are absolutely continuous with respect to
the $d$-dimensional Lebesgue measure $\lbgd$.
Convergence in $\pp(\Rd)$ refers to the weak convergence of measures,
\begin{align}\label{def.weak-convergence}
 \mu_n\to\mu_\infty \quad\mbox{in $\pp(\Rd)$}\quad\text{iff}
 \quad
 \int_X \phi\,\d\mu_n \to \int_X \phi\,\d\mu
 \mbox{ for all $\phi\in C^0_b(\Rd)$.}
\end{align}
Given a measure $\mu\in\pp(\Rd)$, its {\em push-forward}
$\rr_\#\mu\in\pp(X)$ under a Borel map $\rr:\Rd\to \Rd$
is defined by
\begin{align*}
 \rr_\#\mu(A) = \mu\big(\rr^{-1}(A)\big),
\end{align*}
for all Borel sets $A\subset X$.
When $\mu=u\lbgd$ is absolutely continuous with respect to
the Lebesgue measure $\lbgd$ and $\rr$ is an injective
differentiable map with non singular differential $\lbgd$-a.e.,
then $\tilde\mu=\rr_\#\mu\ll\lbgd$ is concentrated on $\rr(\Rd)$
and its density $\tilde\w=\rmd \tilde\mu/\rmd \lbgd$
can be expressed through the formula
\begin{align}
 \label{eq.pushforward}
 \tilde \w(y)=
 \frac{u}{\big| {\rm det} \df \rr\big|}
 \circ \rr^{-1}(y) \quad
 \text{for $\lbgd$-a.e. $y\in \rr(\Rd)$}.
\end{align}
Notice that the particular choice of the dilation map $\rr(x)=rx$, $r>0$,
corresponds to
\eqref{eq:114}, i.e.\ $\rr_\#\mu=(r\id)_\#\mu=(\Dil r u)\lbgd$.

For measures $\mu_0,\,\mu_1\in\pp_2(\setR^d)$,
denote by $W_2(\mu_0,\mu_1)$ their $L^2$-Wasserstein distance
\begin{align}
 \label{eq.transport}
 W_2^2(\mu_0,\mu_1) &= \min_{\gamma\in\Gamma(\mu_0,\mu_1)}
 \int_{\setR^d\times\setR^d} |x-y|^2 \,\d\gamma(x,y),
\end{align}
where $\Gamma(\mu_0,\mu_1)$ is the set of all couplings between $\mu_0$ and $\mu_1$,
i.e.\ all probability measures $\gamma$ on $\setR^d\times\setR^d$
whose marginals are $\mu_0$ and $\mu_1$ respectively.
We collect in the following Theorem some useful properties of
the Wasserstein distance, whose
proof can be found, e.g., in Villani's book~\cite{villani}.
\begin{theorem}
 \label{thm.wasserstein}
 $(\pp_2(\setR^d),W_2)$ is a complete metric space and
 a sequence $\mu_n\in \pp_2(\Rd)$ converges to
 $\mu\in \pp_2(\Rd)$ if and only if
 it converges \eqref{def.weak-convergence} in $\pp(\Rd)$ and
 $\QMom {\mu_n}\to\QMom \mu$.
 $W_2$ is {\em lower semicontinuous} in each argument
 with respect to convergence in $\pp(\setR^d)$,
 and it is {\em convex} in each argument,
 i.e. for $\bar{\mu},\,\mu_0,\,\mu_1\in\pp_2(\setR^d)$ and $\theta\in[0,1]$,
 one has
 \begin{align*}
   W_2\big(\bar{\mu},(1-\theta)\mu_0+\theta\mu_1\big)
   &\leq (1-\theta) W_2(\bar{\mu},\mu_0)
   +\theta W_2(\bar{\mu},\mu_1).
 \end{align*}
 If $\mu_0$ is absolutely continuous with respect $\lbgd$,
 then the minimum in \eqref{eq.transport}
 is achieved by a unique coupling $\gamma_{opt}$
 which is induced by a
 {\em transport map} $\tmap=\tmap_{\mu_0}^{\mu_1}:\setR^d\to\setR^d$,
 i.e.\ $\gamma_{opt}=(\id,\tmap)_\#\mu_0$, $\mu_1=\tmap_\#\mu_0$, and
 \begin{equation}
   \label{eq:32}
   W_2^2(\mu_0,\mu_1)=\int_{\Rd}|\tmap(x)-x|^2\,d\mu_0(x).
 \end{equation}
 In this case the curve $\mu_\vartheta:=\big((1-\vartheta)\id+\vartheta\, \tmap\big)_\# \mu_0$
 for $\vartheta\in [0,1]$ provides a (constant speed, minimizing) geodesic
 connecting $\mu_0$ and $\mu_1$, i.e.\ satisfying
 \begin{equation}
   \label{eq:48}
   W_2(\mu_{\vartheta_0},\mu_{\vartheta_1})=|\vartheta_0-\vartheta_1|
   W_2(\mu_0,\mu_1)\quad
   \text{for all}\quad
   \vartheta_0,\vartheta_1\in [0,1].
 \end{equation}
\end{theorem}
Following McCann \cite{MC}, a functional $\U:\pp_2(\Rd)\to (-\infty,+\infty]$
with proper domain $\dom(\U)\subset \pp_2^r(\Rd)$ is
called \emph{displacement $\kappa$-convex}
(or \emph{geodesically $\kappa$-convex}) if
\begin{equation}
 \label{eq:49}
 \U(\mu_\vartheta)\leq (1-\vartheta)\U(\mu_0)+\vartheta\U(\mu_1)
 -\frac\kappa2 (1-\vartheta) \vartheta W_2^2(\mu_0,\mu_1)\quad
 \forall\, \vartheta\in [0,1],\ \mu_0,\mu_1\in \dom(\U),
\end{equation}
where $\mu_\vartheta$ is a geodesic connecting $\mu_0,\mu_1$ as in \eqref{eq:48}.


\subsection{Integral functionals in $\pp_2(\Rd)$.}
Since the natural framework will be the
space $\pp_2(\Rd)$ of probability measures,
we will use calligraphic letters $\mathcal H$, $\mathcal F$,
to denote the natural extension of
the entropy $\ent$
and the information functional $\fish$ to $\pp_2(\Rd)$,
with the convention that they take the value $+\infty$ when
are evaluated on a measure $\mu\not\in \pp_2^r(\Rd)$.

For every measure $\mu=\w\lbgd\in \pp_2^r(\Rd)$ we thus set
\begin{align}
 \label{eq.energy}
 \Fishm_\al[\mu] := \fish_{\al} (\w)
 = \frac1{2\alpha}
 \int_{\setR^d} \big|\df (\w^{\alpha})\big|^2 \,\dx + \frac\lambda2 \int_{\setR^d} |x|^2 \w(x)\,\dx ,
 \quad
 \text{provided $\w^\alpha\in W^{1,2}(\setR^d)$},
\end{align}
otherwise $\Fishm_\al[\mu]=+\infty$.
Likewise, the associated entropy functional is given by
(recalling our convention that $\Mass u=1$)
\begin{align}
 \label{eq.entropy}
 \rel_\al[\mu] := \ent_{\al} (\w)
 &= \frac{\Theta_\alpha}{\alpha-1/2}
 \Big( \int_{\setR^d} \w^{\alpha+1/2}\,\dx - 1
 \Big)
 + \frac{\Lambda_{\al}}2 \int_{\setR^d} |x|^2 \w \,\dx
 \quad\text{if $\alpha>1/2$} , \\
 \label{eq:50}
 \rel_{1/2,\lambda}[\mu] := \ent_{1/2,\lambda} (\w)&
 = \Theta_{1/2}
 \int_{\setR^d} \w\log \w\,\dx
 + \frac{\Lambda_{1/2,\lambda}}2 \int_{\setR^d} |x|^2 \w \,\dx ,
\end{align}
provided $\w^{\alpha+1/2}\in L^1(\setR^d)$, or $\w\log\w\in L^1(\setR^d)$, respectively,
and $\rel_\al[\mu]=+\infty$ otherwise.
Notice that $\rel_{1/2,0}$ coincides with
the usual logarithmic entropy functional $\frac 12\rel$ associated to
\eqref{eq.logentropy}.
Our choice
$\Lambda_\al=\lambda^{1/2}\big((2\alpha-1)d+2\big)^{-1/2}\geq0$
is motivated by Corollary \ref{cor:ent-fish} and Remark \ref{rem:clarify}.

Some preliminary properties of the information functional $\Fishm_\al$ are collected below.
Since $\mu\in \pp_2(\Rd)$ and
$\Fishm_\al[\mu]=\Fishm_{\alpha,0}[\mu]+\frac \lambda 2 \QMom \mu$,
many properties of $\Fishm_\al$ can be reduced
to the corresponding ones of $\Fishm_{\alpha,0}$.
%

\begin{lemma}[Properties of the perturbed information functionals]
 \label{lem.energy}
 For every $\alpha\in [1/2,1]$ and $\lambda\geq0$
 the functional $\Fishm_\al$ is
 {\em lower semicontinuous} with respect to convergence in $\pp(\setR^d)$,
 and for every sequence $\mu_n=\sigma_n^{1/\alpha}\lbgd$ with $\sigma_n\in W^{1,2}(\Rd)$,
 weakly converging to $\mu$ on $\pp(\Rd)$ with $\sup_n \Fishm_{\alpha,0}[\mu_n]<+\infty$,
 one has
 \begin{equation}
   \label{eq:73}
   \mu=\sigma^{1/\alpha}\lbgd,\quad
   \sigma\in W^{1,2}(\Rd),\quad
   \sigma_n\to \sigma \ \text{strongly}\ \text{in }L^2(\Rd)
   \text{ and weakly in }W^{1,2}(\Rd).
 \end{equation}
 $\Fishm_\al$ is
 {\em convex}\footnote{We recall that $\Fishm_\al$ is {\em not} geodesically convex.}
 (strictly, if $\alpha>1/2$),
 i.e. for $\mu_0,\mu_1\in\pp_2(\setR^d)$ and $\theta\in(0,1)$,
 one has
 \begin{align*}
   \Fishm_\al\big[(1-\theta)\mu_0+\theta\mu_1\big]
   &\leq(1-\theta) \Fishm_{\alpha,\lambda}[\mu_0] +
   \theta \Fishm_{\alpha,\lambda}[\mu_1],
 \end{align*}
 with equality only if $\mu_0=\mu_1$, or if $\alpha=1/2$.
 Moreover, $\Fishm_\al[\mu]$ is finite if and only if
 one of the following two (equivalent) conditions holds:
 \begin{enumerate}
 \item
   $\mu=\ww\lbgd$ with $\ww\in W^{1,1}_{\rm loc}(\Rd)$ and
   $\ww^{\alpha-1}|\df \ww|\in L^2(\Rd)$;
 \item $\mu=\ww\lbgd$ and $\ww^{\alpha+1/2}\in W^{1,1}_{\rm loc}(\Rd)$
   and $\ww^{-1/2}|\df \ww^{\alpha+1/2}|\in L^2(\Rd).$
 \end{enumerate}
 If this is the case,
 then $\ww\in W^{1,p}(\Rd)$ for $p=2/(3-2\alpha)$,
 $\ww^{\alpha+1/2}\in W^{1,1}(\Rd)$, and $\Fishm$ admits the equivalent
 representations
 \begin{equation}
   \label{eq:35}
   \Fishm_{\alpha,0}[\mu]=\frac 1{2\alpha} \int_{\Rd}|\df \ww^\alpha|^2\,\dx=
   \frac \alpha{2}
   \int_{\{\ww>0\}} \ww^{2(\alpha-1)} |\df\ww |^2 \,\dx =
   \Theta_\alpha^2\int_{\{\ww>0\}} \ww^{-1}|\df \ww^{\alpha+1/2}|^2\,\dx.
 \end{equation}
\end{lemma}
\begin{proof}
 We start by proving lower semicontinuity of $\Fishm_{\alpha,0}$.
 Let a sequence $\mu_n=\sigma_n^{1/\alpha}\lbgd$ of measures in $\pp(\Rd)$ with
 limit $\mu$ have
 $\sigma_n\in W^{1,2}(\Rd)$ and $\sup_n \Fishm_{\al}[\mu_n]<+\infty$.
 By Rellich's compactness theorem,
 there exists a subsequence
 (still denoted by $\mu_n$) with $\sigma_n$ strongly converging in
 $L^2_{\rm loc}(\Rd)$ to $\sigma\in L^\alpha(\Rd)$ as $n\to+\infty$.
 It follows that $\mu=\sigma^{1/\alpha}\lbgd$, and
 \begin{displaymath}
   \liminf_{n\to+\infty} \int_\Rd |\df \sigma_n|^2\,\dx\geq
   \int_\Rd |\df \sigma|^2\,\dx.
 \end{displaymath}
 Lower semicontinuity of the perturbed functional
 $\Fishm_\al$ (with $\lambda\geq0$) follows from the lower semi-continuity of the second moment $\QMom \cdot$.

 It follows further that the sequence $\sigma_n$ above
 converges strongly to $\sigma$ in $L^{1/\alpha}(\Rd)$,
 since $\int_\Rd \sigma_n^{1/\alpha}\,\dx = 1$.
 As $1/\alpha\leq 2$, and $\sigma_n$ is bounded in $L^{2^*}(\Rd)$,
 with $2^*=2d/(d-2)> 2$,
 the Sobolev embedding theorem implies strong convergence
 of $\sigma_n$ to $\sigma$ in $L^2(\Rd)$.

 In order to prove the second part of the statement and \eqref{eq:35},
 observe that for $\sigma=\ww^{\alpha}\in W^{1,2}(\Rd)$
 it is easily checked that both $\ww=\sigma^{1/\alpha}$ and $\ww^{\alpha+1/2}=\sigma^{1+1/2\alpha}$
 belong to $W^{1,1}_{\rm loc}(\Rd)$, and that (1)--(2) are satisfied.

 Conversely, if $0<\eps\leq \ww(x)\leq \eps^{-1}$ $\lbgd$-a.e.\ in $\Rd$,
 then the three conditions $\ww\in W^{1,1}_{\rm loc}(\Rd)$,
 $\ww^{\alpha+1/2}\in W^{1,1}_{\rm loc}(\Rd)$,
 and $\ww^{1/\alpha}\in W^{1,1}_{\rm loc}(\Rd)$ are equivalent,
 and any of them implies
 \begin{displaymath}
   \df\ww^{\alpha}=\alpha \ww^{\alpha-1}\df \ww=\frac \alpha{\alpha+1/2}
   \ww^{-1/2}\df \ww^{\alpha+1/2} ,
 \end{displaymath}
 showing \eqref{eq:35} in this case.
 By a standard truncation and approximation argument,
 this property extends to the general situation.

 Concerning the integrability of $\df \ww$,
 we apply H\"older's inequality to the product
 $\df \ww= (\ww^{\alpha-1}\df \ww) \ww^{1-\alpha}$,
 \begin{align*}
   \int_{\Rd} | \df\ww |^p \,\dx
   \leq \Big( \int_{\Rd} (\ww^{\alpha-1}\df\ww)^2\,\dx\Big)^{\frac{p}{2}} \Big( \int_{\Rd} \ww\,\dx \Big)^{\frac{2-p}{2}} .
 \end{align*}
 Next, we observe that $\ww^{\alpha-1}\df\ww\in L^2(\Rd)$ when $\Fishm_\al[\ww]$ is finite,
 and trivially $\ww\in L^1(\Rd)$.
 From here, the Sobolev embedding theorem implies that $\ww\in L^p(\Rd)$.
 An analogous argument, based on
 the splitting $\df \ww^{\alpha+1/2}=
 \ww^{1/2}\big(\ww^{-1/2} \df \ww^{\alpha+1/2}\big)$
 with $\ww^{1/2}\in L^2(\Rd)$ and $\ww^{-1/2}| \df \ww^{\alpha+1/2}| \in L^2(\Rd)$,
 shows that $\ww^{\alpha+1/2}\in W^{1,1}(\Rd)$.

 It remains to prove convexity of $\Fishm_{\alpha,0}$.
 To this end, observe that the function
 \begin{align*}
   (x,\yy)\in\setR_+\times\Rd \, \mapsto \, x^r|\yy|^2
 \end{align*}
 is jointly convex iff $r=2(\alpha-1)\in [-1,0]$;
 this convexity is inherited by $\Fishm_{\alpha,0}$,
 thanks to the second representation in \eqref{eq:35}.
 In fact, the map above possesses a positive definite Hessian,
 and hence is strictly convex, for $-1<r<0$, corresponding to $1/2<\alpha<1$.
 In the borderline case $\alpha=1$,
 the strict convexity of $\int |\df\ww|^2\,\dx$ is obvious.
 The situation $\alpha=1/2$ has been discussed by Gianazza, Savar\'e and Toscani \cite{GST}.
 Finally, $\Fishm_{\alpha,\lambda}$ differs from $\Fishm_{\alpha,0}$ only
 by $\frac{\lambda}{2}\QMom \cdot$, which is convex itself.
\end{proof}
We recall that the Wasserstein slope of a functional
$\mathcal G:\pp_2(\Rd)\to (-\infty,+\infty]$
is defined (as in any metric space) by the formula
\begin{equation}
 \label{eq:37}
 |\partial\mathcal G|(\mu):=
 \limsup_{\nu\to\mu}
 \frac{\big(\mathcal G[\mu]-\mathcal G[\nu]\big)^+}{W_2(\mu,\nu)}\quad
 \forall\,\mu\in \dom(\mathcal G).
\end{equation}
\begin{corollary}[Fundamental entropy--information relation]
 \label{cor:ent-fish}
 For every $\mu\in \pp_2(\Rd)$ with $\Entm_{\alpha,0}[\mu]<+\infty$ we have
 \begin{equation}
   \label{eq:36}
   \Fishm_{\alpha,0}[\mu]=
   \big|\partial\Entm_{\alpha,0}\big|^2(\mu).
 \end{equation}
 Moreover, for every $\lambda\geq0$
 \begin{equation}
   \label{eq:38}
   \Fishm_\al[\mu]=
   \big|\partial\Entm_{\alpha,\lambda}\big|^2(\mu)+
   (\delta_\alpha-2)\Lambda_\al\, \Entm_\al[\mu].
 \end{equation}
\end{corollary}
\begin{proof}
 By \cite[Theorem 10.4.6]{AGS} the slope
 $|\partial\Entm_\al|(\mu)$ is finite if and only if
 $\mu=\ww\lbgd$ with $\ww^{\alpha+1/2}\in W^{1,1}(\Rd)$
 and there exists a $\xxi\in L^2(\mu;\Rd)$ such that
 \begin{equation}
   \label{eq:39}
   \ww\xxi=\Theta_\alpha\df \ww^{\alpha+1/2};\quad
   \text{in this case}\quad
   |\partial\Entm_\al|^2(\mu)=\int_\Rd |\xxi+\Lambda x|^2\,\d\mu.
 \end{equation}
 In the unperturbed case $\lambda=\Lambda=0$,
 relation \eqref{eq:36} follows from
 the third representation of $\Fishm_{\alpha,0}$ given in \eqref{eq:35}.
 In the general case, since $\xxi =
 \Theta_\alpha\ww^{-1/2}\df \ww^{\alpha+1/2}$ with
 $\int_\Rd |\xxi|^2\ww\,\dx=\Fishm_{\alpha,0}[\mu]<+\infty$ and $\delta_\alpha-2=(2\alpha-1)d$,
 we obtain
 \begin{align*}
    \big|\partial\Entm_{\alpha,\lambda}\big|^2(\mu)&=
    \int_{\ww>0}|\xxi+\Lambda x|^2\ww\,\dx=
   \int_\Rd |\xxi|^2\ww\,\dx
   +2\Lambda \Theta_\alpha\int_\Rd x\cdot \df\ww^{\alpha+1/2}\,\dx
   +\Lambda^2 \QMom \mu\\&=
   \Fishm_{\alpha,0}[\mu]-2\Lambda d \Theta_\alpha \int_\Rd \ww^{\alpha+1/2}\,\dx
   +\Lambda^2 \QMom \mu\\&=
   \Fishm_{\alpha,0}[\mu]-\Lambda (\delta_\alpha-2)\Entm_{\al}[\mu]
   +\frac 12{\Lambda^2}\delta_\alpha \QMom \mu\\&=
   \Fishm_{\alpha,\lambda}[\mu]-\Lambda (\delta_\alpha-2)\Entm_{\al}[\mu],
 \end{align*}
 where we applied the relation $\Lambda^2=\lambda/\delta_\alpha$ and the integration by parts formula
 \begin{equation}
   \label{eq:63}
   \int_\Rd x\cdot \df g(x)\,\dx=-d\int_\Rd g(x)\,\dx
   \quad\text{to }g=\ww^{\alpha+1/2}\in W^{1,1}(\Rd).\qedhere
 \end{equation}
\end{proof}
\begin{remark}
 \label{rem:clarify}
 \upshape
 The specific choice of $\Lambda$ in \eqref{eq:38} is probably best understood
 from an alternative derivation of \eqref{eq:38},
 which only uses relation \eqref{eq:36} and the scaling invariance \eqref{eq:115}.
 %
 To avoid technical details,
 we shall sketch this argument in an analogous finite-dimensional framework.
 Consider a smooth function $\V_0:\setR^m\to \setR$ which is \emph{positively homogeneous of degree} $-\delta/2+1$.
 Hence, its squared gradient $\U_0(x)=|\nabla\V_0(x)|^2$ is homogeneous of degree $-\delta$
 (homogeneity in the euclidean setting corresponds to the scaling
 invariance \eqref{eq:115} in the Wasserstein space).
 Setting $\V_\Lambda(x)=\V_0(x)+\frac12\Lambda |x|^2$,
 and observing that $\nabla\V_0(x)\cdot x=(-\delta/2+1)\V_0(x)$ by homogeneity,
 one finds
 \begin{align*}
   |\nabla \V_\Lambda(x)|^2&=|\nabla\V_0(x)+\Lambda x|^2=
   |\nabla\V_0(x)|^2+\Lambda^2 |x|^2+2\Lambda \nabla\V_0(x)\cdot x\\&=
   \U_0(x)+\Lambda^2 |x|^2-\Lambda (\delta-2)\V_0(x)=
   \U_0(x)+\frac12\Lambda^2\delta |x|^2-\Lambda(\delta-2)\V_\Lambda(x).
 \end{align*}
 Therefore, with $\lambda:=\Lambda^2 \delta$ and $\U_\lambda(x):=\U_0(x)+\frac 12\lambda|x|^2$,
 one arrives at the analogue of \eqref{eq:38},
 \begin{displaymath}
   \U_\lambda(x)=|\nabla\V_\Lambda(x)|^2+\Lambda(\delta-2)\V_\Lambda(x).
 \end{displaymath}
\end{remark}

\subsection{A metric characterization of the Heat/Porous medium flow}
For later reference,
we recall the ``Wasserstein characterization'' of the heat/porous medium
equation, 
and its relationship to displacement convexity of the functionals $\rel_{\alpha,\lambda}$.

The following theorem summarizes some of the results from
McCann \cite{MC}, Otto \cite{otto}, Ambrosio, Gigli, and Savar\'e \cite{AS}; see also
Carrillo, McCann, and Villani \cite{CarrilloMcCannVillani06} and Sturm \cite{Sturm05}.
For notational convenience, we introduce the right time derivative
\begin{equation}
 \label{rightd}
 \frac{\rmd^+}{\rmd t}\zeta(t)=
 \limsup_{h\downarrow0}\frac{\zeta(t+h)-\zeta(t)}h
\end{equation}
for real functions $\zeta:[0,+\infty)\to\setR$.
\begin{theorem}
 \label{thm.metricheat}
 The functional $\rel_\al$ is
 displacement $\Lambda$-convex, in the sense of \eqref{eq:49}.
 It generates a $\Lambda$-flow in $\pp_2(\Rd)$:
 for a given initial condition
 $\mu_0=u_0\lbgd\in \dom(\rel_\al)\subset \pp_2^r(\setR^d)$,
 there exists a unique, locally Lipschitz curve
 $\mu:(0,+\infty)\to \pp_2(\Rd)$ with
 $\mu_t=u_t\lbgd\in \dom(\Fishm_\al)\subset \dom(\rel_\al)$
 for every $t>0$,
 which satisfies
 \begin{equation}
   \label{eq:51}
   \frac 12\frac{\d^+}{\d t}W_2^2(\mu_t,\nu)+\frac\Lambda2
   W_2^2(\mu_t,\nu)+\rel_\al[\mu_t]\leq \rel_\al [\nu]
   \quad\forall\, \nu\in \dom(\rel_\al),\ t\geq0,
 \end{equation}
 and attains the initial condition, $\lim_{t\downarrow0}\mu_t=\mu_0$ in $\pp_2(\Rd)$.
 The associated densities $u_t$ constitute
 a weak solution of the porous medium equation
 \begin{align*}
   \partial_t u_t -\Theta\Delta u_t^{\alpha+1/2}+\Lambda\dv(x u_t)=&0
   \quad\text{in }(0,+\infty)\times\Rd,
  \end{align*}
  with the values for $\Theta=\Theta_\alpha>0$ and $\Lambda=\Lambda_\al\geq0$
  given in \eqref{eq.slowdiffparams}.
\end{theorem}
\begin{remark}
 \label{rem:heat}
 \upshape
 The case
 $\alpha=1/2$ and $\lambda=0$ of the previous Theorem concerns
 the relative entropy functional
 \begin{equation}
   \label{eq:56}
   \rel[\mu]=2\rel_{1/2,0}[\mu]:=
   H(u)=\int_\Rd u\log u\,\dx,\quad
   \mu=u\lbgd\in \pp_2^r(\Rd).
 \end{equation}
 Up to a factor two,
 the Wasserstein gradient flow $\mu_t=u_t\lbgd$ of $\rel=2\rel_{1/2,0}$
 corresponds to the heat equation
 \begin{equation}
   \label{eq:52}
   \partial_t u-\Delta u=0\quad\text{in }(0,+\infty)\times\Rd,
 \end{equation}
 whose solution can be expressed by convolution with the initial condition,
 \begin{equation}
   \label{eq:61}
   u_t(x)=\frac 1{(4\pi t)^{d/2}}\int_\Rd \exp(-|x-y|^2/4t)\,\d\mu_0(y).
 \end{equation}
 In particular, $u\in C^\infty((0,+\infty)\times \Rd)$,
 and it is strictly positive.
\end{remark}

\subsection{First variation in the Wasserstein space.}
\label{subsec:first variation}
As we are concerned with the gradient flow of $\Fishm_\al$,
we need to calculate its first variation in the $L^2$-Wasserstein-metric.
Here, we resort to a very direct approach to calculate the latter,
which was introduced by Jordan, Kinderlehrer and Otto \cite{JKO}.
The following notation is needed.
For any smooth, compactly supported vector field $\vcf$ on $\setR^d$,
there exists an associated smooth flow map $\flw^s:\setR^d\to\setR^d$,
satisfying
\begin{equation}
 \label{eq:41}
 \frac{\d}{\d s} \flw^s(x) = \vcf\circ\flw^s(x), \quad \flw^0(x)=x,
\end{equation}
for all $x\in\setR^d$ and all $s\in\setR$.
\begin{lemma}
 \label{lem.subdiff}
 Let a vector field $\vcf\in C^\infty_0(\setR^d)$ and
 a measure $\mu=\w\lbgd\in\pp_2^r(\setR^d)$
 with $\Fishm_{\alpha,\lambda}[\mu]<\infty$ be given.
 Then the map 
 $s\mapsto\Fishm_{\alpha,\lambda}[\flw^s_\#\mu]$
 is differentiable at $s=0$, and
 \begin{align}
   \label{eq.subdiff}
   \frac{\d}{\d s}\Restr{s=0} \Fishm_\al[\flw^s_\#\mu] &=
   \Nolinl\mu\vcf=\Nolin\mu\vcf + \lambda \int_{\setR^d} x\cdot\vcf(x)\,\d\mu(x),
 \end{align}
 where the nonlinear operator $\nolin$ is defined in terms of
 $\s:=u^\alpha$ as
 \begin{equation}
   \label{eq.n}
   \Nolin\mu\vcf = -\frac1{2\alpha} \int_{\setR^d} \Big(
   \alpha \df(\dv\vcf)\cdot \df \s^2
   + 2 \df\vcf\df\s\cdot \df\s
   + (2\alpha-1) (\dv\vcf)|\df\s|^2 \Big)\,\dx.
 \end{equation}
\end{lemma}
\begin{remark}[First variation and weak formulation]
 \label{rem:link weak form}
 \upshape
 Notice that the functional $\Nolin\mu\vcf$ corresponds
 to the expression $N_\alpha(u;\vcf)$ related to the weak formulation of
 \eqref{eq.main} we discussed in Remark \ref{rem:weaker formulation}:
 in fact
 \begin{equation}
   \label{eq:68}
   \Nolin\mu\vcf=N_\alpha(u;\vcf)
   \quad
   \text{if }\mu=u\lbgd\in \pp^r(\Rd).
 \end{equation}
\end{remark}%
\begin{proof}

 It is sufficient to discuss the case $\lambda=0$.
 The evaluation of the first variation for the linear contribution
 can be found in e.g. \cite{JKO}.

 Introduce the pushed-forward measure and its density,
 \begin{align*}
   \mu_s=\big(\flw^s\big)_\#\mu, \quad \w_s = \frac{\d\mu_s}{\d\lbgd},
 \end{align*}
 and write $\s_s=\w_s^{\alpha}$.
 Moreover, for $y\in\setR^d$, define the volume distortion $V_s(y)={\rm det}\big(\df\flw^s(y)\big)>0$.
 According to formula \eqref{eq.pushforward},
 \begin{align}
   \label{eq.flows}
   \w_s\big(\flw^s(y)\big) = \w(y) \cdot \big(V_s(y)\big)^{-1},
   \quad
   \s_s\big(\flw^s(y)\big) = \s(y) \cdot \big(V_s(y)\big)^{-\alpha}.
 \end{align}
 Changing variables $x=\flw^s(y)$ under the integral, we obtain
 \begin{align*}
   \Fishm_{\alpha,0}[\mu_s]
   = \frac1{2\alpha} \int_{\setR^d} \big|\df_x\s_s(x)\big|^2\,\dx
   = \frac1{2\alpha}
   \int_{\setR^d} \big|\big(\df_x\s_s\big)\circ\flw^s(y)\big|^2\,V_s(y)\,dy .
 \end{align*}
 Taking into account that $\flw^0$ is the identity, it follows
 \begin{align}
   \label{eq.subintegral}
   \frac{\d}{\d s}\Restr{s=0} \Fishm_{\alpha,0}[\mu_s]
   = \frac1\alpha
   \int_{\setR^d} \df_x\s(y) \cdot
   \frac{\d}{\d s}\Restr{s=0}\Big[ \big(\df_x\s_s\big)
   \circ\flw^s(y) \, V_s(y)^{1/2} \Big] \,dy .
 \end{align}
 Differentiation of $\s_s$ in \eqref{eq.flows} with respect to $y$ gives
 \begin{align}
   \big(\df_x\s_s\big)\circ(\flw^s(y)\big)\cdot \df_y\flw^s(y)
   = \df_y\s(y)  V_s(y)^{-\alpha} - \alpha \s(y) V_s(y)^{-(\alpha+1)} \cdot \df_yV_s(y).
 \end{align}
 By multiplying both sides by the factor $V_s(y)^{1/2}$ and
 from the right by the matrix inverse of $\df_y\flw^s(y)$,
 one obtains a representation for the expression in the square brackets in \eqref{eq.subintegral},
 \begin{align}
   \label{eq.flowderivative}
   \big(\df_x\s_s\big)\circ\big(\flw^s(y)\big)V_s(y)^{1/2}
   = V_s(y)^{-(\alpha-1/2)} \big( \df_y\s(y) - \alpha\s(y) \df_y(\log V_s(y))\big)\cdot \big(\df_y\flw^s(y)\big)^{-1}.
 \end{align}
 In order to evaluate the $s$-derivative of the right-hand side at $s=0$,
 we make use of the following elementary identities:
 \begin{align}
   \label{eq.erules}
   V_0(y) = 1, \qquad
   \frac{\d}{\d s}\Restr{s=0} V_s(y) = \dv\vcf(y), \qquad
   \frac{\d}{\d s}\Restr{s=0} \big(\df_y\flw^s(y)\big)^{-1} = -\df_y\vcf(y).
 \end{align}
 The first identity is immediate.
 The second follows from the well-known relation between determinant and trace,
 ${\rm det}({\mathbf 1}+sA)=1+s\,{\rm tr}A+o(s)$ for arbitrary square matrices $A$ as $s\to0$.
 The last identity is a consequence of the differentiation rule for inverse matrices,
 $\d(A^{-1})/\d s = -A^{-1}\cdot(\d A/\d s)\cdot A^{-1}$.
 Applying \eqref{eq.erules} to the time-derivative of equation \eqref{eq.flowderivative},
 one obtains
 \begin{align*}
   &\frac{d}{ds}\Restr{s=0}\big[\df_x\s_s\big(\flw^s(y)\big)V_s(y)^{1/2}\big] \\
   &\qquad = -(\alpha-1/2){\dv}_y\vcf(y)\,\df_y\s(y)-\alpha \s(y)\,\df_y{\dv}_y\vcf(y)-\df_y\s(y)\cdot \df_y\vcf(y).
 \end{align*}
 Inserting this into the right-hand side of \eqref{eq.subintegral}
 yields the desired result \eqref{eq.subdiff}.
\end{proof}
Another variation of $\Fishm_\al$ of subsequent relevance is the one
along the dilation group $\flw^s(x)=e^{-s}x$.
A direct computation yields
\begin{lemma}
 \label{le:noncompact}
 The derivative of $\Fishm_\al$ along the dilation group
 $\flw^s(x)=e^{-s}x$ is given by
 \begin{align}
   \label{eq:der dil}
   \frac{\d}{\d s}\Restr{s=0} \Fishm_{\alpha,0}[\flw^s_\#\mu] =
   \scaling_\alpha \Fishm_{\alpha,0}[\mu],\qquad
   \frac{\d}{\d s}\Restr{s=0} \Fishm_\al[\flw^s_\#\mu] =
   \scaling_\alpha \Fishm_{\alpha,0}- \lambda\QMom \mu,
 \end{align}
 where $\scaling_\alpha:=(2\alpha-1)d+2$ is the coefficient
 introduced
 in \eqref{eq.slowdiffparams}.
\end{lemma}
\begin{proof}
 If $\mu=u\lbgd$ with $\sigma=u^\alpha\in W^{1,2}(\Rd)$ then,
 recalling \eqref{eq.pushforward}, $\flw^s_\#\mu=u_s\lbgd$
 with $u_s$ and $\sigma_s=u_s^\alpha$ satisfying
 \begin{displaymath}
   u_s(x)=e^{d\,s}u(e^s x),\quad
   \sigma_s(x)=e^{\alpha d\,s}\sigma(e^sx),\quad
   \df\sigma_s(x)=e^{(\alpha d+1)s}\df\sigma(e^sx),
 \end{displaymath}
 so that a simple computation yields
 \begin{equation}
   \label{eq:85}
   \Fishm_{\alpha,0}[\flw^s_\#\mu]=
   e^{\scaling_\alpha s}\Fishm_{\alpha,0}[\mu],\quad
   \QMom {\flw^s_\#\mu}=e^{-2s}\QMom \mu.
 \end{equation}
 Differentiating \eqref{eq:85} we obtain \eqref{eq:der dil}.
\end{proof}
Notice that \eqref{eq:der dil} is consistent with \eqref{eq.subdiff},
applied to $\zzeta(x)=-x=-\df(|x|^2/2)$,
which generates the dilation flow $\flw^s$.
Although $\zzeta$ is not compactly supported,
Lemma \ref{le:noncompact} above can be derived from Lemma \ref{lem.subdiff}
by an approximation argument.
Fix a smooth real function $\zeta$ such that
$\zeta(r)\equiv 1$ if $r\leq 1$ and $\zeta(r)\equiv0$ if $r\geq 2$ and
we take a family of smooth functions
\begin{equation}
 \label{eq:72}
 \zeta_R(x):=-\frac 12|x|^2 \zeta(|x|/R),\quad
 \zzeta_R(x)=\df \zeta_R(x)= -\bar\zeta(|x|/R)x\quad
 \text{where}\
 \bar\zeta(r)=\zeta(r)+\frac 12 r\zeta'(r),
\end{equation}
so that
\begin{equation}
 \label{eq:79}
 -\zeta_R(x)\uparrow \frac 12|x|^2,\quad
 \zzeta_R(x)\to -x\quad\text{as $R\uparrow+\infty$ $\forall\, x\in \Rd$}.
\end{equation}
It is not difficult to check that
\begin{equation}
 \label{eq:86}
 \lim_{R\uparrow+\infty}\Nolin\mu{\zzeta_R}=\frac
 {\scaling_\alpha}{2\alpha}\int_\Rd|\df \sigma|^2\,\dx=\scaling_\alpha\Fishm_{\alpha,0}[\mu].
\end{equation}

\section{Auxiliary estimates for metric gradient flows}
\label{sec:basic estimate}
This section contains two essential estimates.
The first inequality \eqref{eq:13},
which we will refer to as ``flow interchange estimate'',
gives a precise meaning to the following observation:
{\em the dissipation of one functional along the gradient flow of another functional
 equals the dissipation of the second functional along the gradient flow of the first.}
This estimate will be applied repeatedly in the next section.
The second inequality \eqref{eq:27} provides an equilibration rate
for a gradient flow whose potential is the dissipation
of another functional along its own gradient flow.
Notice that, in view of Corollary \ref{cor:ent-fish},
the information $\Fishm_\al$ and the entropy $\rel_\al$ are connected exactly in this way.

\subsection{Estimates in the smooth finite dimensional case.}

In order to clarify the main ideas,
we shall first derive analogous estimates in a smooth, finite dimensional setting.
Suppose that $\U,\V:\setR^m\to \setR$ are smooth functions and
let us consider the associated gradient flows $S^\U,S^\V:
[0,+\infty)\times \setR^m\to \setR^m$. For every
$u,v\in \setR^m$ the curves $u_t=S^\U_t(u), v_t:=S_t^\V(v)$
are defined as the solutions of the ordinary differential equations
\begin{equation}
 \label{eq:1}
 u_t'=-\nabla \U(u_t),\quad
 v_t'=-\nabla\V(v_t)\quad\text{with initial condition}\quad
 u_0=u,\ v_0=v.
\end{equation}
It is a trivial calculation to check
that $\V$ is decaying along $S^\U$
at precisely the same rate that $\U$ is decaying along $S^\V$.
A quantitative estimate can be obtained in terms of
the rate of dissipation of $\U$ along $S^\V$:
\begin{equation}
 \label{eq:11}
 \rmD^\V\U[v]:=-\frac {\rmd}{\rmd t}\U(S^\V_t(v))\restr{t=0} .
\end{equation}
\begin{lemma}
 The solution $u_t=S^\U_t(u)$ of the $\U$-gradient flow
 satisfies
 \begin{equation}
   \label{eq:3}
   \frac{\rmd }{\rmd t}\V(u_t)+\rmD^\V\U[u_t]=0.
 \end{equation}
 Suppose moreover that $\V$ is $\kappa$-convex, $\rmD^2 \V\geq \kappa\, I$,
 with some constant $\kappa> 0$,
 and that $\U$ can be written as
 \begin{equation}
   \label{eq:2}
   \U(u)=|\nabla\V(u)|^2+2(\theta-\kappa)\, \V(u)
 \end{equation}
 with some $\theta>0$.
 Then $u_t$ converges to the unique minimum point
 $u_{\mathrm{min}}$ of $\V$ (and $\U$) as $t\uparrow+\infty$
 with explicit exponential rates,
 \begin{align}
   \label{eq:4}
   \frac\kappa2|u_t-u_{\mathrm{min}}|^2\leq
   \V(u_t)-\V(u_{\mathrm {min}})&\leq 
   \big(\V(u_0)-\V(u_{\mathrm {min}})\big) \,\exp(-4\kappa\theta\,t),\\
   \label{eq:5}
   \U(u_t)-\U(u_{\mathrm{min}})&\leq 
   \big(\U(u_0)-
 \U(u_{\mathrm{min}})\big)\,\exp\big(-4 (\kappa\land\theta)\theta\,t\big).
 \end{align}
\end{lemma}
\begin{proof}
 \eqref{eq:3} follows immediately by the identity
 \begin{align}\label{eq.bis1}
   \rmD^\V\U[u_t]&= \langle \nabla\U(u_t),
   \nabla\V(u_t)\rangle=
   \lim_{h\downarrow0} h^{-1}\big(\V(u_t)-\V(S^\U_h(u_t))\big)\\&=
   \lim_{h\downarrow0} h^{-1}\big(\V(u_t)-\V(u_{t+h})\big)=
   -\frac \rmd{\rmd t}\V(u_t).
   \label{eq.bis2}
 \end{align}
 Now assume that $\V$ is $\kappa$-convex, and $\U$ is given by \eqref{eq:2}.
 The uniform convexity of $\V$
 entails the existence of a unique minimizer $u_{\mathrm{min}}$
 satisfying
 \begin{equation}
   \label{eq.dan1}
   \frac\kappa2|u-u_{\mathrm{min}}|^2\leq
   \V(u)-\V(u_{\mathrm {min}})\leq \frac 1{2\kappa}|\rmD\V(u)|^2\quad
   \text{for all}\quad
   u\in \setR^m.
 \end{equation}
 Thanks to \eqref{eq:2}, and since $\theta>0$,
 $u_{\mathrm{min}}$ is also the unique minimizer of $\U$.
 A direct computation shows
 \begin{displaymath}
   -\rmD^\V\U(u)=2\big\langle
   \rmD^2 \V(u)\rmD\V(u)+(\theta-\kappa) \rmD\V(u),
   \rmD\V(u)\big\rangle\geq
   2\theta|\rmD\V(u)|^2\geq
   4\kappa\theta\big(\V(u)-\V(u_{\mathrm{min}})\big).
 \end{displaymath}
 An application of Gronwall's estimate to \eqref{eq.bis2} yields \eqref{eq:4}.
 In order to show the second estimate \eqref{eq:5},
 define
 \begin{displaymath}
   \Psi(u):=|\rmD\V(u)|^2,\quad
   \rmD \Psi(u)=2\rmD^2\V(u)\rmD\V(u).
 \end{displaymath}
 Since $\U$ satisfies \eqref{eq:2}, it follows that
 \begin{align*}
   &-\frac{\rmd}{\rmd t}\U(u_t)=|\rmD \U(u_t)|^2=
   |\rmD \Psi(u_t)-2\kappa \,\rmD\V(u_t)+2\theta\,\rmD\V(u_t)|^2
   \\&=
   |\rmD \Psi(u_t)-2\kappa\, \rmD\V(u_t)|^2+4\theta\,\langle \rmD \Psi(u_t)-
   2\kappa\, \rmD \V(u_t),\rmD\V(u_t)\rangle+
   4\theta^2 |\rmD\V(u_t)|^2\geq 4\theta^2 \Psi(u_t),
 \end{align*}
 where we have used that
 \begin{displaymath}
   \langle \rmD \Psi(u_t)-
   2\kappa\, \rmD \V(u_t),\rmD\V(u_t)\rangle=
   2\langle \rmD^2\V(u_t)\rmD\V(u_t),\rmD\V(u_t)\rangle-
   2\kappa|\rmD\V(u_t)|^2\geq0.
 \end{displaymath}
 Since, moreover,
 \begin{equation}
   \label{eq.dan2}
   \U(u)-\U(u_{\mathrm{min}})=\Psi(u)+2(\theta-\kappa)\big(\V(u)-
   \V(u_{\mathrm {min}})\big)\leq
   (1+\kappa^{-1}(\theta-\kappa)^+)\Psi(u)=
   \frac{\theta\lor\kappa}\kappa\Psi(u),
 \end{equation}
 we eventually obtain
 \begin{equation}
   \label{eq:34}
   -\frac{\rmd}{\rmd t}\big(\U(u_t)-\U(u_{\mathrm{min}})\big)\geq
   \frac{4\theta^2\kappa}{\theta\lor\kappa}
   \big(\U(u_t)-\U(u_{\mathrm{min}})\big)=
   4\theta (\kappa\land\theta)\big(\U(u_t)-\U(u_{\mathrm{min}})\big).\qedhere
 \end{equation}
\end{proof}
\begin{remark}
 \label{rem:sharp}
 \upshape
 Estimate \eqref{eq:4} is sharp:
 one can simply take $m=1$ and $\V(u):=\frac \kappa2 u^2$, so that
 $\U(u)=  \kappa \theta u^2$ and $u_t=u_0\exp(-2 \kappa\theta\, t)$.
 The same example shows that also \eqref{eq:5} is sharp, at least in the case
 $\theta\geq \kappa$.
\end{remark}

\subsection{Gradient flows in metric spaces.}
\label{subsec:GF}
We shall now extend the estimates of the previous section to functionals $\U,\V$
and their gradient flows defined in a \emph{complete} metric space $(\rmX,d)$.
Since the previous calculations heavily rely on the use of gradients and scalar products,
it is not immediate how to generalize the results in the absence of a linear structure.

\subsubsection*{\bfseries $\kappa$-convexity and $\kappa$-flows}
Let $\V:\rmX\to (-\infty,+\infty]$ be a proper and lower semicontinuous functional,
with proper domain $\dom(\V)=\{x\in \rmX:\V(x)<+\infty\}$.
We will work under the additional assumption that $\V$ is geodesically $\kappa$-convex:
every $u_0,u_1\in \dom(\V)$ can be connected by
a (minimal, constant speed) geodesic $\vartheta\in [0,1]\mapsto u_\vartheta\in \dom(\V)$
such that
\begin{equation}
 \label{eq:45}
 d(u_{\vartheta},u_{\eta})\leq
 |\vartheta-\eta|d(u_0,u_1),\quad
 \V(u_\vartheta)\leq (1-\vartheta)\V(u_0)+
 \vartheta\V(u_1)-\frac\kappa2\vartheta(1-\vartheta)d^2(u_0,u_1)
\end{equation}
for every $\theta,\eta\in [0,1]$.
Geodesic convexity with $\kappa>0$ implies the existence of
a unique minimizer $u_{\mathrm {min}}$ for $\V$.
In fact, recall the definition of the slope from \eqref{eq:37},
\begin{equation}
 \label{eq:8}
 |\partial \V|(u):=\limsup_{v\to u}
 \frac{\big(\V(u)-\V(v)\big)^+}{d(u,v)}.
\end{equation}
Then $u_{\mathrm {min}}$ satisfies
\begin{equation}
 \label{eq:46}
 \frac \kappa2d^2(u,u_{\mathrm {min}})\leq
 \V(u)-\V(u_{\mathrm {min}})\leq
 \frac 1{2\kappa} |\partial\V|^2(u).
\end{equation}

A gradient flow of $\V$, which is referred to as a $\kappa$-flow,
is a continuous semigroup
$\sgrp\in C^0\big(\dom(\V);\dom(\V)\big)$, i.e.,
\begin{equation}
 \label{eq:6}
 \Sgrp h{\Sgrp tu}=\Sgrp{t+h}u,\quad
 \lim_{h\downarrow0}\Sgrp hu=u\quad\text{for all}\quad
 u\in \dom(\V),\,t\geq0,
\end{equation}
which satisfies the evolution variational inequality
\begin{equation}
 \label{eq:7}
  \frac{1}{2}\frac{\rmd^+}{\rmd t}d^2(\Sgrp tu,v)+
   \frac \kappa 2d^2(\Sgrp tu,v)+\V(\Sgrp tu)
   \leq \V(v) \quad\text{for all } u,v \in \dom(\V),\ t\geq0.
\end{equation}
Equivalently,
the $\kappa$-flow of $\V$ is characterized in Ambrosio and Savar\'e \cite[Thm. 5.7]{AS} by
\begin{equation}
   \label{eq:18}
    |\partial\V|^2(u)=\lim_{h\downarrow0}\frac{\V(u)-
     \V(\Sgrp h u)}h =
   \lim_{h\downarrow0}\frac{d^2(\Sgrp h u,u)}{h^2}
   \quad
   \text{for all}\quad
   u\in \dom(\V).
 \end{equation}

\subsubsection*{\bfseries Minimizing movements.}
We shall now consider another lower semicontinuous functional
$\U: \rmX\to(-\infty,+\infty]$ such that $\dom(\U)\subset \dom(\V)$.
But we will not impose any particular convexity requirement on $\U$.
Thus, the construction of its gradient flow (which is not guaranteed to exist a priori)
is much more difficult in general,
and one cannot expect to find an appealing metric characterization like \eqref{eq:7}.
Therefore, we shall not work on the continuous level;
instead, we reproduce the decay estimates from the finite dimensional situation
at the time-discrete level,
in the framework of so-called \emph{minimizing movements} described in \cite[Chap. 2]{AGS}.

Our main assumption is that
for every time step $\tau\in(0,\tau_o]$,
and every $\bar U\in \rmX$, the functional defined on $\rmX$
\begin{equation}
 \label{eq:57}
 V\mapsto \frac 1{2\tau}d^2(V,\bar U)+\U(V)\quad\text{admits a minimum
   point in $\dom(\U)$.}
\end{equation}
Assumption \eqref{eq:57} expresses that time-discrete gradient flows of $\U$
can be obtained by means of the implicit Euler scheme.
The discrete solution is constructed as follows:
let a partition $\mathcal P_\sttau$ of the time interval $(0,+\infty)$ be given,
with an associated sequence $\sttau=(\tau_n)_{n\in\setN}$ of time steps,
\begin{equation}
 \label{eq:106bis}
 \mathcal P_\sttau:=\{0=t^0_\sttau<t^1_\sttau<\cdots<t^n_\sttau<\cdots\},\quad
 \tau_n=t^n_\sttau-t^{n-1}_\sttau,\quad
 \lim_{n\to+\infty}t^n_\sttau=\sum_n\tau_n=+\infty.
\end{equation}
For later reference, introduce the projection of an arbitrary $t\geq0$ onto the partition by
\begin{align*}
 t_\sttau = \min\{ t^n_\sttau \geq t \,|\, t^n_\sttau \in \mathcal P_\sttau \} .
\end{align*}
Given further an initial datum $U^0_\sttau\in \rmX$,
we can recursively define a sequence $(U^n_\sttau)_{n\in \setN}$ in $\rmX$,
such that $U^n_\sttau$ minimizes the functional \eqref{eq:57}
with $\bar U:=U^{n-1}_\sttau$, i.e.
\begin{equation}
 \label{eq:10}
 \frac 1{2\tau}d^2(U^n_\sttau,U^{n-1}_\sttau)+\U(U^n_\sttau)\leq
 \frac1{2\tau}d^2(V,U^{n-1}_\sttau)+\U(V)
 \quad
 \text{for all}\quad
 V\in \dom(\U).
\end{equation}
Finally, in analogy to \eqref{eq:11} we set
\begin{equation}
 \label{eq:12}
 \rmD^\V\U(u):=\limsup_{h\downarrow0} \frac{\U(u)-\U(\Sgrp hu)}h.
\end{equation}
We are now in the position to formulate and prove the analogous estimates
to \eqref{eq:3}, \eqref{eq:4} and \eqref{eq:5}.
\begin{theorem}
 \label{thm:metric_estimate}
 Let $(U^n_\sttau)_{n\geq 0}$ be an arbitrary solution of \eqref{eq:10}.
 Then, for every $n\geq 1$,
 the following ``flow interchange estimate'' holds:
 \begin{equation}
   \label{eq:13}
   \V(U^n_\sttau)+\tau_n\rmD^\V\U(U^n_\sttau)+
    \frac \kappa2 d^2(U^n_\sttau,U^{n-1}_\sttau)\leq \V(U^{n-1}_\sttau).
 \end{equation}
 Moreover, assume that $\kappa>0$,
 and that
 \begin{equation}
   \label{eq:14}
   \U(u)=|\partial\V|^2(u)+2(\theta-\kappa)\V(u)
   \quad
   \text{for all}\quad
   u\in \dom(\U),
 \end{equation}
 with some $\theta>0$.
 Then $\V$ and $\U$ have the same minimum point $\umin$,
 and for every $n\geq1$, the following estimates hold,
 \begin{align}
   \label{eq:20}
   (1+4\kappa\theta\, \tau_n)\big(\V(U^n_\sttau)
   -\V(\umin)\big)&\leq \V(U^{n-1}_\sttau)-\V(\umin),\\
   \label{eq:27}
   (1+4\theta (\kappa\land\theta)\, \tau_n)\big(\U(U^n_\sttau)-\U(\umin)\big)
   &\leq \U(U^{n-1}_\sttau)-\U(\umin).
 \end{align}
\end{theorem}
Before proving Theorem \ref{thm:metric_estimate},
we briefly comment on its applicability.
Lacking any convexity property of $\U$,
the existence of its gradient flow is not guaranteed \emph{a priori}.
However, the construction by means of \eqref{eq:10} obviously provides
a time-discrete approximation of the sought time-continuous flow.
Denote by $U_\sttau:(0,+\infty)\to \dom(\U)$ the piecewise constant
interpolant of the discrete values on the grid $\mathcal P_\sttau$ \eqref{eq:106}
defined by
\begin{equation}
 \label{eq:58}
 U_\sttau (t):= U^n_\sttau\quad\text{if }t\in (t^{n-1}_\sttau,t^n_\sttau].
\end{equation}
A typical strategy to obtain the gradient flow of $\U$
is to pass to the limit as $\sup_n \tau_n\downarrow0$,
in (a suitable subsequence of) the time-discrete flows $U_\sttau$.
In order to obtain the limit by a compactness argument,
$\ttau$-independent \emph{a priori} estimates for the $U_\sttau$ are needed.
Formula \eqref{eq:13} provides a powerful tool to derive such estimates,
as will be seen in section \ref{sct.existence},
where \eqref{eq:13} is applied with $\U=\Fishm_\al$,
and various choices of $\V$
(see also \cite[chapters 2,3]{AGS} for further examples).

On the other hand, estimates \eqref{eq:20} and \eqref{eq:27} are
used to give a quantitative description of the equilibration behavior
of solutions to \eqref{eq.main} in section \ref{sct.asymptotics}.
Notice that by Corollary \ref{cor:ent-fish},
the functionals $\Fishm_\al$ and $\rel_\al$ are related exactly in the way \eqref{eq:14}.


\begin{proof}
 In order to prove \eqref{eq:13},
 choose $V:=\Sgrp h{U^n_\sttau}$ with $h>0$
 in the variational inequality \eqref{eq:10} for $\U$.
 This gives
 \begin{equation}
   \label{eq:15}
   \U(U^n_\sttau)-\U(\Sgrp h{U^n_\sttau})\leq
   \frac 1{2\tau_n}\Big(d^2(\Sgrp h{U^n_\sttau},U^{n-1}_\sttau)-
   d^2(U^n_\sttau,U^{n-1}_\sttau)\Big).
 \end{equation}
 Dividing by $h>0$ and passing to the limit as $h\downarrow0$,
 the variational characterization \eqref{eq:7} of $\V$,
 with $t=0$, $u=U^n_\sttau$ and $v=U^{n-1}_\sttau$, yields
 \begin{equation}
   \label{eq:16}
   \rmD^\V\U(U^n_\sttau)\leq \tau_n^{-1}
       \Big(\V(U^{n-1}_\sttau)-\V(U^n_\sttau)-
       \frac \kappa2 d^2(U^n_\sttau,U^{n-1}_\sttau)\Big),
 \end{equation}
 which is estimate \eqref{eq:13}.

 Let us now suppose that $\U$ is as in \eqref{eq:14}, with $\kappa>0$.
 The existence of a (unique) minimum point $\umin$ of $\V$
 follows easily from the convexity assumption \eqref{eq:45}
 and the completeness of $\rmX$
 (see e.g.\ \cite[Lemma 2.4.8]{AGS}).
 Inequality \eqref{eq:46} and the fact that $|\partial\V|(\umin)=0$
 show that $\umin$ is also the unique minimum point of $\U$.

 From now on, we assume without loss of generality that $\V(\umin)=\U(\umin)=0$.
 Since, by \cite[Thm. 5.7]{AS},
 \begin{equation}
   \label{eq:17}
   \limsup_{h\downarrow0}\frac{|\partial\V|^2(u)-
     |\partial\V|^2(\Sgrp h u)}h\geq 2\kappa |\partial\V|^2(u)
   \quad
   \text{for all}\quad
   u\in \dom(|\partial\V|),
 \end{equation}
 the second characterization \eqref{eq:18} of the $\kappa$-flow
 induces a lower bound on $\rmD^\V\U$,
 \begin{equation}
   \label{eq:19}
   \rmD^\V\U(u)\geq 2\kappa |\partial\V|^2(u)+
   2(\theta-\kappa)|\partial\V|^2(u)=2\theta |\partial\V(u)|^2
   \geq 4\theta\kappa \V(u) , 
 \end{equation}
 where the last inequality follows from the $\kappa$-convexity of $\V$, see \eqref{eq:46}.
 Insert \eqref{eq:19} in \eqref{eq:13} and neglect the
 non-negative term $\kappa d^2(U^n_\sttau,U^{n-1}_\sttau)/2$ to obtain estimate \eqref{eq:20}.

 It remains to prove the last inequality \eqref{eq:27}.
 Estimate the difference on the right hand side of \eqref{eq:15} by the triangle inequality
 to find
 \begin{equation}
   \label{eq:23}
   \U(U^n_\sttau) - \U(\Sgrp h{U^n_\sttau})
   \leq d(\Sgrp h{U^n_\sttau},U^{n}_\sttau)\Big(
   d(\Sgrp h{U^n_\sttau},U^{n-1}_\sttau)+
   d(U^n_\sttau,U^{n-1}_\sttau)\Big).
 \end{equation}
 Divide this inequality by $h$ and pass to the limit $h\downarrow0$,
 then use the second variational characterization \eqref{eq:18} of $\sgrp$
 to identify the limit on the right hand side.
 In combination with \eqref{eq:19}, we obtain the two-sided estimate
 \begin{equation}
   \label{eq:21}
   2\theta\tau_n |\partial\V|^2(U^n_\sttau)\leq
   \tau_n \rmD^\V\U(U^n_\sttau)\leq |\partial\V|(U^n_\sttau)
   d(U^n_\sttau,U^{n-1}_\sttau).
 \end{equation}
 This leads further to
 \begin{equation}
   \label{eq:24}
   2\theta\tau_n |\partial\V|(U^n_\sttau)\leq
   d(U^n_\sttau,U^{n-1}_\sttau), \quad
   4\theta(\theta\land \kappa)\tau_n^2\U(U^n_\sttau)\leq d^2(U^n_\sttau,U^{n-1}_\sttau).
 \end{equation}
 For the derivation of the second estimate above,
 we have made use of (recall $\U(\umin)=0$)
 \begin{equation}
   \label{eq:26}
   \U(u)\leq \frac {\theta\lor \kappa}\kappa |\partial\V|^2(u) ,
 \end{equation}
 which follows in the same way as \eqref{eq.dan1} in the finite-dimensional case,
 simply using \eqref{eq:14} instead of \eqref{eq:2},
 and \eqref{eq:46} instead of \eqref{eq.dan2}.

 Direct substitution of \eqref{eq:24} into the characterization \eqref{eq:10}
 of the time-discrete flow for $\U$, with $V:=U^{n-1}_\sttau$,
 gives a bound similar to \eqref{eq:27} but with
 the smaller constant $2\theta(\theta\land \kappa)$.
 A more refined estimate is needed to recover the better constant.

 To this end,
 introduce the so-called De Giorgi's variational interpolants:
 for $r\in (0,\tau_n]$ we define $U^{n,r}_\sttau$ a minimizer of the functional
 $V\mapsto \frac 1{2 r}d^2(V,U^{n-1}_\sttau)+\U(V)$, thus obtaining
 a path joining $U^{n-1}_\sttau$ with $U^n_\sttau$
 (which can be chosen equal to $U^{n,\tau_n}_\sttau$).
 Replacing $\tau_n$ with $r$, and $U^n_\sttau$ with $U^{n,r}_\sttau$,
 estimates \eqref{eq:24}
 and monotonicity of $\U$, i.e. $\U(U^{n,r}_\sttau)\geq \U(U^n_\sttau)$, see \cite[Lemma 3.1.2]{AGS},
 yield
 \begin{equation}
   \label{eq:47}
    4\theta(\theta\land\kappa)r^2\U(U^{n}_\sttau)
    \leq 4\theta(\theta\land\kappa)r^2\U(U^{n,r}_\sttau)
    \topref{eq:24}\leq
   d^2(U^{n,r}_\sttau,U^{n-1}_\sttau)\quad
   \forall\, r\in (0,\tau_n].
 \end{equation}
 The crucial inequality satisfied by
 $U^{n,r}_\sttau$ is provided by
 \cite[Thm. 3.1.4]{AGS}, which shows that
 the map $r\mapsto d(U^{n,r}_\sttau,U^{n-1}_\sttau)$ is Borel and
 \begin{equation}
   \label{eq:42}
   \U(U^n_\sttau)+\frac 1{2\tau_n}d^2(U^n_\sttau,U^{n-1}_\sttau)+
   \int_0^{\tau_n} \frac 1{2r^2}d^2(U^{n,r}_\sttau,U^{n-1}_\sttau)\,\d r\leq
   \U(U^{n-1}_\sttau).
 \end{equation}
 Inserting \eqref{eq:24} and \eqref{eq:47} into \eqref{eq:42}
 we eventually get \eqref{eq:27}.
\end{proof}
The following remark shows that the discrete estimates
\eqref{eq:20} and \eqref{eq:27}
give rise to equilibration estimates analogous to
\eqref{eq:4} and \eqref{eq:5}.
\begin{remark}
 \label{rem:exp}
 \upshape
 Let $(A^n_\sttau)_{n\geq0}$ be a sequence of non-negative numbers,
 and let $A_\sttau$ be its piecewise constant interpolant,
 taking the value $A^n_\sttau$ in the interval $(t^{n-1}_\sttau,t^n_\sttau]$.
 For a given time $t\geq0$ we set $\sft_\sttau:=\inf\{t^k_\tau:t^k_\tau\geq t, k\in \setN\}$.
 Assume that
 \begin{equation}
   \label{eq:29}
   (1+\mathsf c\,\tau_n)
   A^n_\sttau\leq A^{n-1}_\sttau\quad\text{for all $n\in \setN$ and
     some $\mathsf c\geq0$}.
 \end{equation}
 Then, with $\mathsf c_\sttau(t):=\tau_n^{-1}\log(1+\mathsf c \tau_n)$ for $t\in (t^{n-1}_\sttau,t^n_\sttau]$,
 we get
 \begin{equation}
   \label{eq:30}
   A_\sttau(t)\leq A_\sttau(s)\exp\Big(-\int_{s_\ssttau}^{t_\ssttau} \mathsf c_\sttau(r)\,\d r\Big)
   \quad
   \forall\, t>s>0.
 \end{equation}
 If $A^0_\sttau$ is uniformly bounded,
 Helly's Theorem shows that any  sequence $\ttau_k$
 of partitions admits a subsequence (still denoted by $\ttau_k$)
 and a limit function $A$ such that
 $\lim_{k\uparrow\infty}A_{\sttau_k}(t)=A(t)$ for every $t\geq0$.
 If $\lim_{k\uparrow\infty}t_{\sttau_k}=t$ for every $t\in [0,+\infty)$ then
 we eventually get
 \begin{equation}
   \label{eq:31}
   A(t)\leq A(s)\exp\big(-\mathsf c (t-s)\big).
 \end{equation}
\end{remark}

\section{Existence of Solutions}
\label{sct.existence}

The existence proof follows the ideas first developed
by Jordan, Kinderlehrer and Otto \cite{JKO}
and later generalized in Ambrosio, Gigli and Savar\'e \cite{AGS}.
First, a family of time-discrete approximative solutions
is constructed by means of the variational \emph{minimizing movement} scheme,
which has been introduced in section \ref{subsec:GF} above.
The decay of the information $\Fishm_\al$ is immediate from the construction and leads to a priori estimates.
The latter provide weak compactness of the set of discrete solutions,
and thus allows to conclude the existence of a time-continuous limit curve.

In order to show that the limit curve constitutes a weak solution to the gradient flow,
we need to pass to the limit in the non-linear discrete equation.
Additional compactness is needed for this step.
A sufficient \emph{a priori} estimate is obtained by evaluating the dissipation
of the logarithmic entropy along the discrete flow.

The procedure is thus similar to the one developed in Gianazza, Savar\'e and Toscani \cite{GST}.
However, we will follow a different strategy to derive the \emph{a priori} estimate
and to prove the convergence of the scheme.

\subsection{The semi-discrete scheme}
\label{sct.semidiscrete}
For the reader's convenience,
we review the main steps of the semi-discretization procedure in a series of lemmas.
For any further details,
we refer to the exhaustive treatment in \cite{AGS}.
Given a partition $\mathcal P_\sttau$ of $[0,+\infty)$
induced by the sequence of time-steps $\ttau=(\tau_n)_{n\in\setN}$
as in \eqref{eq:106bis}, and an initial measure $\Mu^0_\sttau
\in\pp_2(\setR^d)$,
we consider the
sequence $(\Mu^n_\sttau)_{n\in \setN}$
recursively defined by solving
the following variational problem in $\pp_2(\setR^d)$:
\begin{equation}
 \label{eq:59}
 \text{find}\quad
 \Mu^n_\sttau\in \pp_2(\Rd)\qquad
 \text{which minimizes}\qquad
 \Mu\mapsto \frac1{2\tau_n}W^2_2(\Mu^{n-1}_\sttau,\Mu) +
 \Fishm_\al[\Mu] .
\end{equation}
The existence and uniqueness of a minimizer follows
by standard methods from the calculus of variations,
employing the continuity and convexity properties of $W_2$ and $\Fishm_\al$
collected in Theorem \ref{thm.wasserstein} and Lemma \ref{lem.energy}.
\begin{lemma}[Basic discrete estimates]
 For each choice of $\ttau$ and $\Mu^0_\sttau\in\pp_2(\setR^d)$,
 there exists a unique sequence
 $(\Mu^n_\sttau)_{n\in \setN}\subset \pp_2^r(\Rd)$ solving \eqref{eq:59}.
 The information functionals $\Fishm_\al[\Mu_\sttau^n]$ are 
 finite and monotone non-increasing with respect to $n$,
 \begin{align}
   \label{eq.apriori1}
   \Fishm_\al[\Mu_\sttau^n]
   \leq \Fishm_\al[\Mu_\sttau^{n-1}] \leq \Fishm_\al[\Mu^0_\sttau].
 \end{align}
 Moreover, the quadratic moments and
 Wasserstein distances between
 measures at consecutive time steps satisfy for every $N>0$
 \begin{align}
   \label{eq.apriori2}
   \sum_{n=1}^\infty
   \tau_n\Big(\frac{W_2(\Mu_\sttau^n,\Mu_\sttau^{n-1})}{\tau_n}
   \Big)^2 \leq 2 \Fishm_\al[\Mu^0_\sttau],
   \qquad
   \sup_{1\leq n\leq N}\QMom {M^n_\sttau}\leq 2\,
   \QMom {M^0_\sttau}+ 4 t^N_\sttau\,\Fishm_\al[M^0_\sttau].
 \end{align}
\end{lemma}
\begin{proof}
 By the minimality property of $\Mu_\sttau^n$, one has
 \begin{align}
   \label{eq.wasserspeed}
   \frac1{2\tau_n}W_2^2(\Mu_\sttau^{n-1},\Mu_\sttau^n) + \Fishm_\al[\Mu_\sttau^n] &\leq \Fishm_\al[\Mu_\sttau^{n-1}],
 \end{align}
 which induces \eqref{eq.apriori1} and the first inequality of
 \eqref{eq.apriori2} by summation over $n=1,2,\ldots$.
 In order to prove the second bound,
 we sum up the triangular inequalities
 \begin{displaymath}
   \sqrt{\Mom2{M^n_\sttau}}\leq
   \sqrt{\Mom2{M^{n-1}_\sttau}}+W_2(M^n_\sttau,M^{n-1}_\sttau),
 \end{displaymath}
 from $n=1$ to $n=N$, observing that
 \begin{displaymath}
   \sum_{n=1}^N W_2(M^n_\sttau,M^{n-1}_\sttau)\leq
   \Big(\sum_{n=1}^N \tau_n\Big)^{1/2}\,\Big(\sum_{n=1}^N
   \tau_n^{-1}W_2^2(M^n_\sttau,M^{n-1}_\sttau)\Big)^{1/2}.
   \qedhere
 \end{displaymath}
\end{proof}
The traditional way to proceed from here is to introduce the discrete velocity vector fields.
However, for our existence proof, we shall not pursue that line of argument.
Nonetheless, we briefly recall the definition and basic properties of the velocity field
for the sake of completeness, and also
since it provides a natural interpretation of the gradient flow structure.

As each $\Mu_\sttau^n$ is absolutely continuous,
Theorem \ref{thm.wasserstein} guarantees the existence
and the ($\Mu^n_\sttau$-essential) uniqueness of a Borel map
$\tmap_\sttau^n:\setR^d\to\setR^d$ with $(\tmap_\sttau^n)_\#\Mu_\sttau^n=\Mu_\sttau^{n-1}$,
realizing the optimal transport from $\Mu_\sttau^n$ to $\Mu_\sttau^{n-1}$
in the definition of the Wasserstein distance \eqref{eq.transport}.
Introduce the {\em discrete velocity vector field} ${\mathbf v}_\sttau^n$ by
\begin{align}
 \label{eq.velo}
 \velo_\sttau^n(x) := \frac1{\tau_n} \big(x-\tmap_\sttau^n(x)\big).
\end{align}
For each $n=1,2,\ldots$,
the $\velo_\sttau^n$ satisfies the nonlinear characterization
\begin{align}
 \label{eq.dnonlinear0}
 \int_\Rd \vcf(x)\cdot\velo_\sttau^n(x)\,\d\Mu_\sttau^n(x) &=
 \Nolinl{\Mu_\sttau^n}\vcf
\end{align}
for all test vector fields $\vcf\in C^\infty_0(\setR^d,\setR^d)$,
with $\nolinl$ as defined in \eqref{eq.n}.
As showed by \cite{JKO},
the relation \eqref{eq.dnonlinear0} follows by taking the first variation
of the minimizing functional \eqref{eq:59}
at the minimum point $\Mu^n_\sttau$ along the flow $\mathsf S_s(\mu):=\flw^s_\#\mu$,
obtained by taking the push-forward of a measure $\mu$
under the maps \eqref{eq:41} induced by the vector field $\vcf$.

The discrete solution thus satisfies a system composed of
the non-linear constraint \eqref{eq.dnonlinear0},
and a transport equation,
which follows directly from the definition \eqref{eq.velo}
and reads in weak form as
\begin{align}
 \label{eq.dtransport}
 \int_\Rd \phi(x)\d \big(\Mu_\sttau^n(x) - \Mu_\sttau^{n-1}(x) \big)
 = \int_\Rd \big( \phi(x) - \phi(x-\tau_n\velo_\sttau^n(x)) \big)\,\d\Mu_\sttau^n(x) .
\end{align}
The system \eqref{eq.dnonlinear0}--\eqref{eq.dtransport}
is a canonical starting point for studying the continuous time limit.

However, we shall take another approach
which does not make use of the velocity vector fields.
Instead, we will derive the weak form of the gradient flow equations
in the continuous time limit by means of the next lemma.
\begin{lemma}[Discrete time derivative]
 \label{lem.velo}
 Let a test function $\zeta\in C^\infty_0(\Rd)$ be given,
 which satisfies $-|\kappa| I\leq D^2\zeta\leq |\kappa| I$ for some $\kappa<0$.
 Then
 \begin{equation}
   \label{eq:40}
   -\frac{|\kappa|}2W_2^2(\Mu^n_\sttau,\Mu^{n-1}_\sttau)\leq
   \int_\Rd \zeta\,\d\Mu^n_\sttau-
   \int_\Rd \zeta\,\d\Mu^{n-1}_\sttau
   + \tau_n\Nolinl{\Mu_\sttau^n}{\df \zeta}
   \leq
   \frac{|\kappa|}2 W_2^2(\Mu^n_\sttau,\Mu^{n-1}_\sttau).
 \end{equation}
 Moreover, the second moment satisfies
 \begin{equation}
   \label{eq:84}
   (1+2\lambda\tau_n)\QMom {M^n_\sttau}+W_2^2(\Mu^n_\sttau,\Mu^{n-1}_\sttau)=
   \QMom {M^{n-1}_\sttau}+2\scaling_\alpha\tau_n
   \Fishm_{\alpha,0}[M^n_\sttau].
 \end{equation}
\end{lemma}
Notice that \eqref{eq:40} is a discrete local version
of the weak formulation \eqref{eq:66bis} for the gradient flow of $\Fishm_\al$.
This can be seen after division of both sides by $\tau_n>0$,
and taking into account that we expect $W_2^2(M_\sttau^n,M_\sttau^{n-1})=o(\tau_n)$
(at least in an integral sense) in the limit $\tau_n\downarrow0$.
\begin{proof}
 Choosing $\vcf:=-\df \zeta$,
 the semi-group $\mathsf S_s^\V$ coincides with
 the $\kappa$-Wasserstein gradient flow of the functional
 $\V(\mu):=\int_\Rd \zeta\,\d\mu$ \cite[Example 11.2.2]{AGS}.
 We apply the ``flow interchange estimate'' \eqref{eq:13} to $\U:=\Fishm_\al$,
 where
 \begin{displaymath}
   \df^\V\U(\mu)=-\Nolinl\mu{\zzeta}=\Nolinl\mu{\df\zeta}
 \end{displaymath}
 by Lemma \ref{lem.subdiff}.
 Substituting $\V$ with $-\V$ yields \eqref{eq:40}.
 The particular choices $\zeta:= \frac 12 \kappa |x|^2$ with $\kappa = \pm 1$
 correspond to $\V(\mu)= \pm\frac 12 \QMom \mu$  in Theorem \ref{thm:metric_estimate}
 and yield \eqref{eq:84} via Lemma \ref{le:noncompact}.
\end{proof}

\subsection{An a priori estimate related to entropy dissipation}
\label{sct.apriori}
The classical estimates \eqref{eq.apriori1}--\eqref{eq.apriori2} on the discrete scheme
are sufficient to prove the existence of continuous limits of $\Mu_\sttau$ and $\velo_\sttau$
as $\sup |\tau_n|\to0$.
An additional a priori estimate is needed to verify that these limits
indeed satisfy the desired nonlinear evolution equation \eqref{eq:66bis}.
Below, we obtain a sufficient estimate
from the decay of the logarithmic entropy $\rel = 2 \rel_{1/2,0}$ along the discrete flow.

To simplify notations, we introduce
\begin{align}
 \label{eq.beta}
 \beta:=\frac{1-\alpha}\alpha,
\end{align}
which is in one-to-one correspondence to $\alpha$,
and varies between $\beta=0$ for $\alpha=1$, and $\beta=1$ for $\alpha=1/2$.
\begin{theorem}
 \label{thm:crucialapriori}
 Let $\Mu^n_\sttau=U^n_\sttau\,\lbgd$ be any solution
 of the discrete variational scheme \eqref{eq:59},
 and define $\sigma^n_\sttau=(U^n_\sttau)^{\alpha}$.
 Then $U^n_\sttau\in W^{2,1}(\Rd)$ and $\sigma^n_\sttau\in W^{2,2}(\Rd)$
 satisfy for every $n\geq 1$
 \begin{align}
   \label{eq:crucialapriori}
   \sfc_0\tau_n
   \int_{\setR^d} \big\| \df^2 \sigma^n_\sttau \big\|^2 \,\dx
   \leq  \rel [\Mu^{n-1}_\sttau] - \rel[\Mu^{n}_\sttau] +
   d\lambda \tau_n,\quad\text{where }
   \sfc_0 := \frac1\alpha\Big( 1 - \frac{(d-1)^2}{d(d+2)}\beta \Big).
 \end{align}
\end{theorem}
\begin{proof}
 The expression on the left-hand side of \eqref{eq:crucialapriori} constitutes a lower bound
 on the dissipation rate of $\rel$ along the gradient flow of $\Fishm_\al$.
 In the spirit of Theorem \ref{thm:metric_estimate},
 the roles of $\rel$ and $\Fishm_\al$ will be interchanged:
 we calculate the dissipation of the information $\Fishm_\al$ along the
 Wasserstein gradient flow of $\rel$, which corresponds to the
 classical heat equation (see Remark \ref{rem:heat}).
 The regularizing effect of the heat flow
 makes it comparatively easy to justify the necessary manipulations.

 Apply the ``flow interchange estimate'' \eqref{eq:13}
 with the choice $\V:=\rel$ (which is geodesically convex, thus $\kappa=0$)
 and $\U:=\Fishm_\al$.
 Then \eqref{eq:crucialapriori} follows once we have shown that
 \begin{equation}
   \label{eq:60}
   \df^\V \U(\mu)\geq \sfc_0
   \int_{\setR^d} \big\| \df^2 \sigma \big\|^2 \,\dx- d\lambda \quad
   \text{if }\mu=u\lbgd,\ \sigma=u^\alpha.
 \end{equation}
 This is the content of Lemma \ref{le:technical} below.
\end{proof}
A remark concerning the choice of the logarithmic entropy $\rel$ 
as a Lyapunov functional is due at this point.
First, it is a convenient choice since the resulting entropy dissipation
has the right homogeneity to provide $W^{1,2}$-compactness for $\sigma$ easily.
Second, it is a canonical choice since the heat flow generated by $\rel$
dissipates every functional $\Psi:\pp(\Rd)\to\setR$ that is jointly convex
in $u$ and $\df u$; so in particular it dissipates $\Fishm_\al$ for $1/2\leq\alpha\leq1$.
\begin{lemma}
 \label{le:technical}
 Let $\mu_0=u_0\lbgd \in \dom(\Fishm)$ and let $\mu_s=u_s\lbgd$
 be the associated solution to the heat equation \eqref{eq:52} on $\Rd$,
 which is given by \eqref{eq:61}.
 If the right derivative
 $\liminf_{h\downarrow0} \frac1h\big(\Fishm_\al[\mu] -
 \Fishm_\al[\mu_h]\big)$ is finite, then $\sigma_0=u_0^\alpha\in W^{2,2}(\Rd)$ and
 \begin{equation}
   \label{eq.alongheat}
   -\frac {\d^+}{\d s}\Fishm_\al[\mu_s]
   \restr{s=0}=\liminf_{h\downarrow0} \frac1h\big(\Fishm_\al[\mu_0] - \Fishm_\al[\mu_h]\big) \geq
   \sfc_0  \int_\Rd \big\| \df^2 \sigma_0 \big\|^2 \,\dx
   -d\lambda.
 \end{equation}
\end{lemma}
\begin{remark}
 \label{rem:easy_case}
 \upshape
 We emphasize that the dimension-dependent prefactor of $\beta$ in the definition of $\sfc_0$
 is always less than one, and converges to one for $d\to\infty$.
 Hence, for $0\leq\beta\leq1$ --- corresponding to $\frac12\leq\alpha\leq1$ --- the coefficient $\sfc_0$ is positive.
 Further, observe that in the case $\alpha=1$, it follows $\sigma=u$,
 and \eqref{eq.alongheat} reduces for $\lambda=0$ to the well known estimate
 \begin{displaymath}
   -\frac {\d^+}{\d s}\frac 12\int_\Rd |\df u_s|^2\,\d x=
   \int_\Rd |\Delta u_0|^2\,\dx=
   \int_\Rd |\df^2 u_0|^2\,\dx.
 \end{displaymath}
\end{remark}
\begin{proof}
 Since
 $\Fishm_\al[\mu]=\Fishm_{\alpha,0}[\mu]+\frac\lambda2\QMom \mu$
 and
 \begin{equation}
   \label{eq:62}
   \frac {\d}{\d s}\QMom {\mu_s}=
   \frac {\d}{\d s} \int_\Rd |x|^2 \,\d\mu_s =
   -2\int_\Rd x\cdot \df u_s\,\dx\topref{eq:63} =
   2 d\int_\Rd u_s\,\dx=2 d,
 \end{equation}
 thanks to $u_s\in W^{1,1}(\Rd)$,
 we can simply consider the case of $\lambda=0$
 and evaluate the time derivative
 of $\Fishm_{\alpha,0}[\mu_s]=\frac 1{2\alpha}\int_\Rd |\df u_s^\alpha|^2\,\dx$.

 By standard parabolic theory,
 $u_s$ is a $C^\infty$-smooth,
 strictly positive probability density for each $s>0$.
 Moreover, for every $\bar s>0$ there exists a constant $C_{\bar s}$ such
 that for all $s\geq\bar s$ and  $p\in [1,+\infty]$,
 \renewcommand{\ww}{u}
 \begin{equation}
   \label{eq:43}
   \begin{aligned}
     \|\ww_s\|_{L^p(\Rd)}+ \|\df\ww_s\|_{L^p(\Rd;\Rd)}+
     \|\df^2 \ww_s\|_{L^p(\Rd;\setR^{d\times d})}&\leq C_{\bar s}.
   \end{aligned}
 \end{equation}
 We stress that this classical estimate holds 
 with any probability density $\ww_0$ as initial condition,
 thanks to the representation \eqref{eq:61} of solutions to the heat equation.
 For $\eps>0$ we consider the smooth real functions
 \begin{equation}
   \label{eq:25}
   f_\eps(r):=(\eps^{1/\alpha}+r)^\alpha-\eps,\quad r\in [0,+\infty),
 \end{equation}
 whose first and second derivatives are uniformly bounded;
 recall that $1/2\leq \alpha\leq1$.
 Then the functions $\sw_{\eps,s}:=f_\eps(\ww_s)\in C^\infty(\Rd)$
 satisfy bounds analogous to \eqref{eq:43}
 (with constants also depending on $\eps>0$)
 and the evolution equation
 \begin{align}
   \label{eq.heatskeweps}
   \partial_s \sw_{\eps} = \Delta \sw_\eps
   + \beta (\eps+\sw_\eps)^{-1} | \df\sw_\eps |^2 ,
 \end{align}
 where $\beta$ has been defined in \eqref{eq.beta}.
 It follows that for every test function $\zeta\in C^\infty_0(\Rd)$,
 \begin{align*}
   &- \frac 1{2\alpha}\frac{\d}{\d s}\int_\Rd |\df \sw_{\eps,s}|^2\zeta\,\dx=
   - \frac1\alpha \int_\Rd \df\sw_{\eps}\cdot\partial_s\df\sw_{\eps,s}\, \zeta\,\dx \\
   &= \frac1\alpha \int_{\Rd} \Delta\sw_\eps \Big( \Delta \sw_\eps +
   \beta (\eps+\sw_\eps)^{-1}|\df\sw_\eps|^2 \Big)\,\zeta\,\dx
   +\frac1\alpha\int_\Rd \Big( \Delta \sw_\eps +
   \beta (\eps+\sw_\eps)^{-1}|\df\sw_\eps|^2 \Big)\df\sw_\eps\cdot \df\zeta\,\dx
 \end{align*}
 We integrate the previous inequality with respect to
 $s$ from arbitrary points
 $0<s_1<s_2$ and we choose a family of cutoff functions
 $\zeta(x)=\zeta_n(x):=\zeta_0(x2^{-n})$, where $0\leq \zeta_0\leq 1$ and $\zeta_0(x)=1$ if $|x|\leq 1$.
 The contribution of last integral vanishes in the limit as $n\to+\infty$,
 since $\df \zeta_n(x)=2^{-n} \df \zeta_0(x2^{-n})$.
 Consequently,
 \begin{align}
   \label{eq:44}
   -\frac 1{2\alpha}\frac{\d}{\d s}\int_\Rd |\df \sw_{\eps,s}|^2\,\d x
   &= \frac1\alpha \int_{\Rd} \Big(\big(\Delta\sw_\eps\big)^2
   +\beta\Delta \sw_\eps
    (\eps+\sw_\eps)^{-1}|\df\sw_\eps|^2 \Big)\,\dx\\
    \label{eq:55}&=
    \frac1\alpha \int_{\Rd} \Big(\big(\Delta\sw_\eps\big)^2
    +4\beta\Delta \sw_\eps
    |\df z_\eps|^2 \Big)\,\dx.
 \end{align}
 where $z_\eps:=(\eps+\sw_\eps)^{1/2}-\eps^{1/2}$
 satisfies $2\df z_\eps=(\eps+\sw_\eps)^{-1/2}\df \sw_\eps$.


 For further estimation, we recall that
 \begin{equation}
   \label{eq:53}
   0 =
   \int_\Rd \big\|\df^2\sw_\eps\big\|^2\,\dx -
   \int_\Rd \big(\Delta\sw_\eps\big)^2\,\dx,
 \end{equation}
 and the following integration by parts rule
 (see \cite[Theorem 3.1]{GST} for a proof),
 \begin{align}
   \label{eq:54}
   0 &= \int_{\Rd} \dv\big(|\df z_\eps|^2 \df \sw_\eps\big)
   \,\dx =
   \int_{\Rd} \Big(
   \Delta\sw_\eps \big|\df z_\eps\big|^2
   + 2 \df^2\sw_\eps \df z_\eps\cdot \df z_\eps
   -  4 \big|\df z_\eps\big|^4
   \Big)\,\dx,
 \end{align}
 Introducing the parameters
 \begin{align}
   \sfc_1 = - \frac{4d\beta}{\alpha(d+2)} = \frac1\alpha\Big(\frac{8\beta}{d+2} - 4\beta\Big), \quad
   \sfc_2 = \sfc_0 + \frac{(d-1)\beta}{\alpha(d+2)} = \frac1\alpha\Big(1+\frac{(d-1)\beta}{d(d+2)}\Big),
 \end{align}
 we add $\sfc_1$ times the equality \eqref{eq:54} and $\sfc_2$ times the equality \eqref{eq:53}
 to the dissipation relation \eqref{eq:55}.
 This gives
 \begin{align}
   \nonumber
   -\frac 1{2\alpha}\frac{\d}{\d s}\int_\Rd |\df \sw_\eps|^2\,\d x=&
   \sfc_0 \int_{\Rd} \|\df^2\sw_\eps\|^2\,\dx
   +  \frac\beta{\alpha(d+2)} \int_{\Rd} \bigg\{
   (d-1)\big\|\df^2\sw_\eps\big\|^2
   - \frac{d-1}{d} \big({\Delta\sw_\eps}\big)^2 \\
   \label{eq.afteribp}
   & - 8d  \,{\df^2\sw_\eps}\df z_\eps\cdot\df z_\eps
   + 8 \,{\Delta\sw_\eps} \big|{\df z_\eps}\big|^2
   + 16 d \,\big|{\df z_\eps}\big|^4
   \bigg\}\,\dx .
 \end{align}
 As a final ingredient, 
 we need the following elementary estimate
 for the error term in the Cauchy-Schwarz inequality
 (see e.g. \cite[Lemma 2.1]{JM} for a proof).
 \begin{lemma}
   \label{lem.matrix}
   Let $\sfA=(a_{ij})$ be a real symmetric $d\times d$-matrix,
   and let $0\neq\ee\in\setR^d$ be an arbitrary vector.
   Then the square-norm $\|\sfA\|^2=\sum a_{ij}^2$ and
   the trace ${\rm tr}\sfA=\sum a_{ii}$ of $\sfA$ satisfy
   \begin{align}
     \label{eqn.matrix}
     (d-1) \|\sfA\|^2 - \frac{d-1}{d} \big({\rm tr}\sfA\big)^2
     \geq \frac1{d\, |\ee|^4} \big( d\, \sfA\,\ee\cdot\ee - {\rm tr} \sfA\,|\ee|^2 \big)^2.
   \end{align}
 \end{lemma}
 Lemma \ref{lem.matrix} is applied with $\sfA=\df^2\sw_\eps$ and $\ee=\df z_\eps$
 to estimate the norm of the Hessian inside the second integral in \eqref{eq.afteribp}.
 A straight-forward calculation reveals
 \begin{align}
   \label{eqn.dissipfinal}
  -\frac 1{2\alpha}\frac{\d}{\d s}\int_\Rd |\df \sw_\eps|^2\,\d x
  \geq \sfc_0 & \int_\Rd \|\df^2\sw_\eps\|^2\,\dx \\
  & + \frac{d\beta}{\alpha(d+2)}
  \int_\Rd \Big(\frac{\df^2\sw_\eps \df z_\eps \cdot \df z_\eps}{|\df z_\eps|^2}
  - \frac1d \Delta\sw_\eps - 4 | \df z_\eps |^2 \Big)^2 \,\dx.
 \end{align}
 In view of $\beta\geq0$, the last term is non-negative.
 Integration of \eqref{eqn.dissipfinal} from $s=0$ to $s=h>0$ thus provides
 \begin{align*}
   \frac1{2\alpha} \int_\Rd |\df \sw_{\eps,h}|^2\,\d x
   + \sfc_0 \int_0^h \int_\Rd \|\df^2\sw_{\eps,s}\|^2\,\dx\,\d s 
   \leq \frac1{2\alpha} \int_\Rd |\df \sw_{\eps,0} |^2\,\d x
 \end{align*}
 By definition of $\sw_{\eps}$ from $\ww$ 
 by means of $f_\epsilon$ given in \eqref{eq:25}, 
 it follows,
 \begin{align*}
   \sw_{\eps,s} = \big( \sigma_s^{1/\alpha} + \eps^{1/\alpha} \big)^\alpha , \quad \sigma_s = \ww_s^\alpha .
 \end{align*}
 Using that $\alpha\leq1$, we conclude by definition of $\Fishm_{\alpha,0}$,
 \begin{align*}
   \frac1{2\alpha} \int_{\Rd} | \df \sw_{\eps,0} |^2 \,\d x 
   = \frac\alpha2 \int_{\Rd} \frac{| \df \sigma_0 |^2}{(\sigma_0^{1/\alpha} + \eps^{1/\alpha})^{2(1-\alpha)}}\,\d x 
   \leq \Fishm_{\alpha,0}[\mu_0].
 \end{align*}
 Moreover, it clearly follows that $\ww_{\eps,s}\lbgd\to\mu_s=\sw_s\lbgd$ 
 in the sense of $\pp(\Rd)$ as $\eps\downarrow0$, pointwise in $s\in[0,h]$.
 Hence, by lower semi-continuity of $\Fishm_{\alpha,0}$ and of the $W^{2,2}$-semi-norm,
 \begin{align*}
   \Fishm_{\alpha,0}[\mu_h]
   + \sfc_0 \int_0^h \int_\Rd \|\df^2\sw_h\|^2\,\d x\,\d s 
   \leq \Fishm_{\alpha,0}[\mu_0] .
 \end{align*} 
 After addition of the contribution due to the confinement potential,
 and division by $h>0$,
 the limit $h\downarrow0$ provides the desired estimate \eqref{eq.alongheat}.
\end{proof}

\subsection{Passage to the continuous time limit}
\label{sct.ctslimit}
To conclude the construction of the discrete approximation,
we introduce interpolants of $\Mu_\sttau^n$ and $U^n_\sttau$ for all times $t\geq0$:
denote by $\Mu_\sttau$, $U_\sttau$, and $\sft_\sttau$, respectively,
the right-continuous piecewise constant functions with
\begin{align*}
 \Mu_\sttau(t) := \Mu_\sttau^n=U^n_\sttau\lbgd, \quad
 U_\sttau(t) := U_\sttau^n,\quad
 \sft_\sttau:=t^n_\sttau
 \qquad\text{for $t\in\big(t^{n-1}_\sttau,t^n_\sttau\big]$.}
\end{align*}
Moreover, introduce accordingly $\sigma_\sttau=(U_\sttau)^\alpha$.

We recall the notion of \emph{generalized minimizing movement}
introduced by {De Giorgi} \cite{DeG} (and here adapted to the case
of a scheme with variable time steps):
\begin{definition}[Generalized minimizing movements]
 \label{def:GMM}
 For a given partition $\mathcal P_\sttau$,
 the function $\Mu_\sttau:[0,+\infty)\to\pp_2(\Rd)$ obtained by piecewise constant interpolation
 of the solutions $\Mu_\sttau^n$ of the minimization scheme \eqref{eq:59}
 are called {\em discrete solutions}.

 We say that a curve $\mu\in C^0([0,+\infty);\pp_2(\Rd))$ is a
 \emph{generalized minimizing movement} for $\Fishm_\al$ with
 initial datum $\mu_0\in \pp_2(\Rd)$
 if there exist
 \begin{enumerate}[(i)]
 \item a sequence of partitions $\ttau_k$ with $\lim_{k\to\infty}\sup\big\{\tau_n:\sum_{i=1}^n\tau_i\leq T\big\}=0$ for every $T>0$;\\
 \item sequence of initial data $M_{\sttau_k}^0$ converging to
   $\mu_0$ in
   $\pp_2(\Rd)$ with
   $\rel_\al[M^0_{\sttau_k}]\to \rel_\al[\mu_0]$ and
   $\Fishm_\al[M^0_{\sttau_k}]\to \Fishm_\al[\mu_0]$;\\
 \item a corresponding sequence of discrete solutions $M_{\sttau_k}(t)$
   such that
   \begin{equation}
     \label{eq:107}
     M_{\sttau_k}(t)\to \mu(t)\quad\text{in }\pp_2(\Rd)\quad \forall\, t\geq 0.
   \end{equation}
 \end{enumerate}
 We denote by  $GMM(\Fishm_\al;\mu_0)$ the collection of all the generalized minimizing movements.
\end{definition}
Concerning assumption (ii), we remark that 
convergence of $\Fishm_\al$ implies convergence of $\rel_\al$ if $\mu_0\in\dom(\Fishm_\al)$.

The next theorem is the main statement of this section
and in particular yields Theorem \ref{thm.existence} in the case
$\mu_0\in \dom(\Fishm_\al)$.
\begin{theorem}
 \label{thm:convergence}
 Assume that  $\mu_0\in \dom(\Fishm_\al)$.
 For every sequence of partitions $\ttau_k$ and discrete solutions
 $M_{\sttau_k}(t)$ whose initial data $M_{\sttau_k}(0)=M^0_{\sttau_k}$
 satisfy conditions (i) and (ii) of the previous definition,
 there exists a subsequence (still denoted by $\ttau_k$)
 such that as $k\to\infty$,
 \begin{gather}
   \label{eq:69}
   M_{\sttau_k}=U_{\sttau_k}\lbgd\to \mu=u\lbgd\ \text{in
   $\pp_2(\Rd)$, uniformly in each compact interval $[0,T]$} ,\\
   \label{eq:71bis}
   U_{\sttau_k}(t)\to u(t)\quad\text{in }L^{\alpha+1/2}(\Rd)\cap L^1(\Rd)\quad \forall\, t\geq0,\\
   \label{eq:71}
   \sigma_{\sttau_k}\to \sigma\ \text{strongly in
   }L^2((0,T);W^{1,2}(\Rd)),
   \ \text{and weakly in }L^2((0,T);W^{2,2}(\Rd)),
   \ \forall\, T>0,\\
   \label{eq:118}\Fishm_{\al}[M_{\sttau_k}(t)]\to \Fishm_{\al}[\mu(t)]\quad
   \text{for $\mathcal L^1$-a.e.\ $t>0$}.
 \end{gather}
 The limit $u$ is a solution of \eqref{eq.main} in the
 weak form \eqref{eq.weakform}, \eqref{eq:65}.
 In particular, any generalized minimizing movement
 $\mu=u\lbgd\in GMM(\Fishm_\al;\mu_0)$ is a weak solution of \eqref{eq.main} in this sense.
\end{theorem}
\begin{proof}
The proof is naturally divided into various steps.

\noindent
\underline{Step 1: Compactness.}
\emph{We prove that there exists a suitable subsequence (still denoted by $\ttau_k$)
 such that for every $t,T\geq0$}
 \begin{gather}
   \label{eq:69bis}
   M_{\sttau_k}(t)=U_{\sttau_k}(t)\lbgd\to \mu(t)=u(t)\lbgd\quad \text{in
   }\pp(\Rd),\quad
   \sigma_{\sttau_k}(t)\to \sigma(t)\quad\text{strongly in }L^2(\Rd),\\
   \label{eq:70}
   \sigma_{\sttau_k}(t)\to \sigma(t)\quad\text{weakly in }W^{1,2}(\Rd),\quad
   \sigma_{\sttau_k}\to \sigma\quad\text{weakly in }L^2((0,T);
   W^{2,2}(\Rd)).
 \end{gather}
Pointwise weak convergence in $\pp(\Rd)$ follows immediately from
\eqref{eq.apriori2} and a simple extension of
Ascoli-Arzel\`a compactness Theorem (see e.g. \cite[Prop. 3.3.1]{AGS}),
which in particular implies
\begin{align}
 \label{eq.holder}
 W_2\big(\Mu_\sttau(t),\Mu_\sttau(s)\big) \leq
 \big(2\Fishm_\al[\Mu^0_\sttau]
 \big)^{1/2}\cdot(\sft_\sttau-s_\sttau)^{1/2}\quad
 \forall\, 0\leq s\leq s_\sttau\leq t\leq \sft_\sttau.
\end{align}
Strong (resp.\ weak) convergence in $L^2(\Rd)$ (resp.\ in $W^{1,2}(\Rd)$)
then follows by the properties of $\Fishm_\al$ stated in Lemma \ref{lem.energy},
and the uniform bound $\Fishm_{\al}[M_\sttau^n]\leq\Fishm_\al[M_\sttau^0]$, see \eqref{eq.apriori1}.

An explicit bound of the quadratic moment, even in the case
$\lambda=0$, can be easily obtained from \eqref{eq:84},
which gives
\begin{equation}
 \label{eq:87}
 \QMom {M_\sttau(t)}\leq \QMom {\mu_0}+2\scaling_\alpha \sft_\sttau \Fishm_{\alpha,0}(M^0_\sttau).
\end{equation}
Finally, weak convergence in $L^2((0,T);W^{2,2}(\Rd))$
is a consequence of the uniform bound
\begin{equation}
 \label{eq:74}
 \sfc_0\int_0^T \int_\Rd \|\df^2 \sigma_\sttau\|^2\,\dx\,\d t\leq
 \rel[M^0_\sttau]+\pi\QMom {M^0_\sttau}+ (d\lambda+2\pi\scaling_\alpha\Fishm_{\alpha,0}[M^0_\sttau]) (T+\tau)\quad
 \forall\, T>0.
\end{equation}
Here and below, $\tau := \sup \{\tau_i \mid i \in \setN\}$ controls the time step size
in the given partition,
and inequality \eqref{eq:74} is obtained by
summing the central \emph{a priori} estimate \eqref{eq:crucialapriori} for $\sigma_\sttau$
with the bound \eqref{eq:84} on the second moment,
taking into account (see e.g.\ \cite[Sect. 2.3]{GST}) that
\begin{equation}
 \label{eq:88}
 \rel[\mu]+\pi\QMom \mu\geq0\quad\forall\, \mu\in \pp_2(\Rd).
\end{equation}

\noindent
\underline{Step 2: Discrete equation.}
\emph{We prove that for every $0\leq s\leq {\mathsf s}_\sttau=s^m_\sttau\leq t\leq \sft_\sttau=t^p_\sttau$
 and $\zeta\in C^\infty_0(\Rd)$ with $-\kappa I\leq \df^2\zeta\leq \kappa I$
 one has}
\begin{equation}
 \label{eq:75}
 \Big|\int_\Rd \zeta\,\d M_\sttau(t)-
 \int_\Rd \zeta \,\d M_\sttau(s)+
 \int_{s_\ssttau}^{t_\ssttau}\Nolinl{M_\sttau(r)}{\df \zeta}\,\d r\Big|\leq
 \kappa \tau \Fishm_\al[M^0_\sttau]
\end{equation}
Summing up \eqref{eq:40} from $n=m+1$ to $n=p$
one gets
\begin{displaymath}
 \Big|\int_\Rd \zeta\,\d M_\sttau^p-
 \int_\Rd \zeta \,\d M_\sttau^m+
 \sum_{n=m+1}^p  \tau_n\, \Nolinl{M_\sttau(r)}{\df \zeta} \Big|\leq
 \frac \kappa2\sum_{n=m+1}^p W^2_2(M^n_\sttau,M^{n-1}_\sttau).
\end{displaymath}
Claim \eqref{eq:75} above then follows recalling \eqref{eq.apriori2}.

\noindent
\underline{Step 3: Continuous limit equation.}
\emph{We prove for the limit curve $\mu$ of step 1,
 that for every $0\leq s\leq t$ and $\zeta\in C^\infty_0(\Rd)$,
 one has
 \begin{equation}
   \label{eq:75bis}
   \int_\Rd \zeta\,\d \mu(t)-
   \int_\Rd \zeta \,\d \mu(s)+
   \int_{s}^{t}\Nolinl{\mu(r)}{\df \zeta}\,\d r=0,
 \end{equation}
 and that the Lebesgue density $u(t,\cdot)$ of $\mu(t)$ satisfies
 the weak formulation \eqref{eq:66bis} and \eqref{eq:65}
 of the gradient flow \eqref{eq.main}.}\\
Equation \eqref{eq:75bis} above follows by passing to the limit in \eqref{eq:75}
as $k\uparrow\infty$ and applying Step 1.
Notice that weak convergence of $\sigma_{\sttau_k}$ to $\sigma$ in $L^2((0,T);W^{2,2}(\Rd))$
and strong convergence in $L^2((0,T);L^2(\Rd))$ immediately yields
strong convergence in $L^2((0,T);W^{1,2}(\Rd))$ and,
up to the extraction of a further subsequence, \eqref{eq:118}.
Further, observe that the map $\sigma\mapsto N(\sigma^{1/\alpha},\df\zeta)$ is continuous and quadratic
in $W^{1,2}(\Rd)$.
Hence, one can pass to the limit in the integral term
$\int_{s_\ssttau}^{t_\ssttau}\Nolinl{M_\sttau(r)}{\df \zeta}\,\d r$.

It is then a standard fact that \eqref{eq:75bis} implies that $\sigma$
satisfies the weak formulation \eqref{eq:66bis}.
Finally, the other weak formulation \eqref{eq:65} can be deduced
integrating by parts, which is admissible since $\sigma=u^\alpha\in L^2_{\mathrm{loc}}([0,+\infty);W^{2,2}(\Rd))$;
see Remark \ref{rem:weaker formulation}.

\noindent
\underline{Step 4: Convergence of quadratic moments.}
\emph{We prove that for every $t\geq0$,
 \begin{equation}
   \label{eq:89}
   \lim_{k\uparrow+\infty}\QMom{M_{\sttau_k}(t)}=\QMom{\mu(t)}\quad
   \text{so that}\quad M_{\sttau_k}(t)\to \mu(t)\quad\text{in }\pp_2(\Rd).
 \end{equation}}%
Insert the family of test functions $\zeta_R$ --- introduced in \eqref{eq:72} ---
in formula \eqref{eq:75bis} above.
Recalling that $\mu(t)\in \pp_2(\Rd)$ and arguing like in the derivation of \eqref{eq:86},
we get
\begin{equation}
 \label{eq:90}
 \QMom{\mu(T)}+2\lambda \int_0^T\QMom{\mu(r)}\,\d r=\QMom{\mu_0}+
 2\scaling_\alpha\int_0^T \Fishm_{\alpha,0}[\mu(r)]\,\d r
\end{equation}
On the other hand, summing up estimate \eqref{eq:84} from $n=1$ to $N$ yields
\begin{equation}
 \label{eq:92bis}
 \QMom{M_\sttau(T)}+2\lambda \int_0^T \QMom{M_\sttau(r)}\,\d r\leq
 \QMom{\mu_0}+2\scaling_\alpha \int_0^{T_\ssttau} \Fishm_{\alpha,0}[M_\sttau(r)]\,\d r
\end{equation}
The strong convergence of $\sigma_{\sttau_k}$ to $\sigma$ in $L^2_{\mathrm{loc}}((0,+\infty);W^{1,2}(\Rd))$
permits passage to the limit
in the right-hand side of \eqref{eq:92bis} as $k\uparrow\infty$.
Convergence of $M_{\sttau_k}(t)$ in $\pp(\Rd)$,
and the lower semicontinuity of the quadratic moment yield for all $T>0$,
\begin{align*}
 \limsup_{k\uparrow+\infty}\QMom{M_{\sttau_k}(T)}
 \leq \QMom{\mu(T)} +  2\lambda \bigg( \int_0^T\QMom{\mu(r)}\,\d r - \liminf_{k\uparrow+\infty}\int_0^T \QMom{M_{\sttau_k}(r)}\,\d r \bigg)
 \leq \QMom{\mu(T)} .
\end{align*}
This implies pointwise convergence
(and a posteriori also locally uniform convergence, see \eqref{eq.holder})
of $M_{\sttau_k}(t)$ in $\pp_2(\Rd)$.

\noindent
\underline{Step 5: Convergence in $L^{\alpha+1/2}(\Rd)$.}
\emph{We prove that the densities $U_\sttau$ converge strongly in $L^{\alpha+1/2}$.}\\
In view of the previously established convergence results,
is obviously sufficient to show that the norms $\|U_\sttau(t)\|_{L^{\alpha+1/2}(\Rd)}$
converge to their respective limits $\|u(t)\|_{L^{\alpha+1/2}(\Rd)}$.
Observe that these norms and the (unperturbed) entropies are equivalent,
\begin{align*}
 \|U_\sttau(t)\|_{L^{\alpha+1/2}}^{\alpha+1/2} = \frac{\alpha-1/2}{\Theta_\alpha} \rel_{\alpha,0}[\Mu_\sttau(t)] +1 .
\end{align*}
By the convergence of $\Mu_{\sttau_k}$ in $\pp_2(\Rd)$
and the uniform boundedness of $\Fishm_{\alpha,0}$,
\begin{equation}
 \label{eq:108}
 \Big|\rel_{\alpha,0}[\Mu_{\sttau_k}(t)]- \rel_{\alpha,0}[\mu(t)]\Big|
 \leq W_2(M_{\sttau_k}(t),\mu(t))\big(\Fishm_{\alpha,0}(\Mu_{\sttau_k}(t))+\Fishm_{\alpha,0}(\mu(t))\big)^{1/2}
\end{equation}
converges to zero, and so does the difference of the respective $L^{\alpha+1/2}(\Rd)$-norms.

Thus, the proof of Theorem \ref{thm:convergence} is complete.
\end{proof}
To complete the proof of Theorem \ref{thm.existence},
it remains to treat the case in which the initial condition $u_0$
has finite logarithmic entropy $\ent$, but infinite information.
\begin{theorem}[Regularizing effect]
 Suppose $\mu_0 = u_0 \lbgd \notin \dom (\Fishm_\al)$ but
 \eqref{eq:83} is satisfied, and
 the discrete initial data $M^0_\sttau=U^0_\sttau\lbgd$ converging to $\mu_0$ in
 $\pp_2(\Rd)$ satisfy $\rel[M^0_{\sttau_k}]\to \ent[u_0] < +\infty$.
 Then the conclusions of Theorem \ref{thm:convergence} remain valid
 even in the absence of hypothesis (ii) of Definition \ref{def:GMM},
 except that \eqref{eq:69} and \eqref{eq:71bis} may fail at the endpoint $t=0$.
\end{theorem}
\begin{proof}
 A combination of Sobolev inequalities and standard interpolation results yields
 \begin{equation}
   \label{eq:88bis}
   \forall\, \eps>0\ \exists\,Q_\eps>0:\quad
   \|\df \sigma\|^2_{L^2(\Rd)}\leq \eps\, \|\df^2 \sigma\|_{L^2(\Rd)}^2+
   Q_\eps\,\|\sigma\|_{L^{1/\alpha}(\Rd)}^2
 \end{equation}
 for every function $\sigma\in  W^{2,2}(\Rd)\cap L^{1/\alpha}(\Rd)$; recall that $1/\alpha\leq2$.

 Multiply the bound \eqref{eq:84} on moments by $2\pi$,
 and take the sum with the a priori estimate \eqref{eq:crucialapriori},
 \begin{align*}
   \rel[M^n_\sttau]+2\pi\QMom{M^n_\sttau}+\sfc_0\tau_n
   \|\df^2\sigma^n_\sttau\|_{L^2(\Rd)}^2&\leq
   \rel[M^{n-1}_\sttau]+2\pi\QMom{M^{n-1}_\sttau}\\
   &\quad+\tau_n\Big(\lambda d
   +
   2\alpha^{-1}\pi\scaling_\alpha \|\df\sigma^n_\sttau\|_{L^2(\Rd)}^2\Big) .
 \end{align*}
 Choosing $\eps>0$ in estimate \eqref{eq:88bis} such that
 $2\alpha^{-1}\pi\scaling_\alpha \eps= \sfc_0/2$,
 and observing that $\|\sigma^n_\sttau\|_{L^{1/\alpha}(\Rd)}=1$,
 we obtain
 \begin{equation}
   \label{eq:90bis}
   \rel[M^n_\sttau]+2\pi\QMom{M^n_\sttau}+\frac{\sfc_0}2\tau_n
   \|\df^2\sigma^n_\sttau\|_{L^2(\Rd)}^2\leq
   \rel[M^{n-1}_\sttau]+2\pi\QMom{M^{n-1}_\sttau}
   +\sfC\tau_n,
 \end{equation}
 where $\sfC:=\lambda d+ 2 \alpha^{-1} \pi \scaling_\alpha Q_\eps$.
 In particular, in view of inequality \eqref{eq:88} this implies
 the estimate
 \begin{equation}
   \label{eq:91}
   \frac{\sfc_0}2\int_0^T\big\|\df^2 \sigma_\sttau(t)\|_{L^2(\Rd)}^2\,\d t
   \leq
   \rel[M^0_\sttau]+2\pi\QMom{M^0_\sttau}+\sfC T_\sttau.
 \end{equation}
 Now combine the fact that $n\mapsto \Fishm_\al [M^n_\sttau]$ is nonincreasing
 with the moment estimate \eqref{eq:92bis} and inequality \eqref{eq:91} above
 to deduce both a pointwise and an integral bound for
 the functional $\Fishm_\al$:
 \begin{align*}
   T\,\Fishm_\al[\Mu_\sttau(T)]
   & \leq \int_0^T \Fishm_\al[\Mu_\sttau(t)]\,\d t \\
   & \leq (1+\delta_\alpha/2)\int_0^{T_\sttau} \Fishm_{\alpha,0}[\Mu_\sttau(t)]\,\d t + \frac14\QMom{\Mu_\sttau^0} \\
   & \leq \frac{1+\delta_\alpha/2}{2\alpha} \int_0^{T_\sttau} \big( \eps \|\df^2\sigma_\sttau(t)\|^2_{L^2(\Rd)} + Q_\eps \big)\,\d t
   + \frac14 \QMom{\Mu_\sttau^0} \\
   & \leq \sfC_1\big(\rel[M^0_\sttau]+2\pi\QMom{M^0_\sttau} \big) +\sfC_2 \,(1+T_\sttau).
 \end{align*}
 With these estimates at hand,
 we can repeat the arguments given in the proof of Theorem \ref{thm:convergence}
 on each interval $[\eps_h,+\infty)$,
 for a vanishing sequence $\eps_h>0$.
 A standard diagonal argument concludes the proof.
\end{proof}

\section{Convergence to Equilibrium}
\label{sct.asymptotics}

The large time behavior of the previously constructed weak solutions $\w(t)$
to \eqref{eq.main} is described as follows.
For a positive confinement force $\lambda>0$,
the solutions converge to a Barenblatt profile $\blatt_\al$,
which is the common minimizer of the entropy $\ent_\al$ and the information $\fish_\al$,
in $L^1(\Rd)$ at the exponential rate of $\lambda$;
see section \ref{sct.confined} below.
In addition to the convergence estimates in $L^1(\Rd)$,
we obtain estimates for the equilibration of the information $\fish_\al$.
For $\alpha=1$,
this translates into a convergence estimate for $\w(t)$ in $W^{1,2}(\Rd)$;
see section \ref{sct.w12}.

For vanishing confinement,
no minimizer of $\fish_{\alpha,0}$ exists, and solutions converge weakly to zero.
However, the solutions become asymptotically self-similar:
they converge in $L^1(\Rd)$ after a suitable rescaling.
The limit is again a Barenblatt profile $\blatt_\al$.
The rescaling arguments are spelled out in section \ref{sct.intermediate}.

To close the discussion, we provide an alternative (but largely formal) argument
that leads to $L^1(\Rd)$ convergence in section \ref{sct.ct}.
The argument extends the idea developed by Carrillo and Toscani \cite{CT2}
for the confined thin film equation ($\alpha=1$ and $\lambda=1$) in $d=1$ dimension.

\subsection{Equilibration of the confined gradient flow}
\label{sct.confined}
Throughout this section, we assume positive confinement $\lambda>0$.

In order to study the equilibration of solutions to \eqref{eq.main},
we shall employ the estimates \eqref{eq:20} and \eqref{eq:30}
that have been derived in the general metric framework in section \ref{sec:basic estimate}.
At the very basis of our argument is the intimate relation \eqref{eq.specialform}
between the entropy $\rel_\al$ and the information $\Fishm_\al$.
We recall their definitions,
\begin{align*}
 \rel_\al[\mu] := \ent_{\alpha,\lambda}[\w]
 = \frac {\Theta_\alpha}{\alpha-1/2}\bigg( \int_{\setR^d} \w^{\alpha+1/2}\,\dx - \Big(\int_{\setR^d} \w\,\dx\Big)^{\alpha+1/2} \bigg)
 + \frac{\Lambda_\al}{2} \int_{\setR^d} |x|^2\,\mu(\dx),
\end{align*}
where
\begin{align}
 \label{eq.Lambda}
 \Theta_\alpha =\frac{\sqrt{2\alpha}}{2\alpha+1},\qquad
 \Lambda_\al = \sqrt{\frac{\lambda}{\delta_\alpha}},\qquad
 \delta_\alpha=(2\alpha-1)d+2.
\end{align}
And, according to Corollary \ref{cor:ent-fish},
the information can be expressed in terms of the entropy as follows,
\begin{equation}
 \label{eq:38bis}
 \Fishm_\al[\mu]=
 \big|\partial\Entm_{\alpha,\lambda}\big|^2(\mu)+
 (\delta_\alpha-2)\Lambda\, \Entm_\al[\mu].
\end{equation}
Thus, $\Fishm_\al$ is connected to $\rel_\al$ in precisely the same way
that $\U$ is connected to $\V$ in the equation \eqref{eq:14}
considered in section \ref{sec:basic estimate};
see \eqref{eq:77} below for the exact correspondences.

Since $\lambda>0$, also $\Lambda=\Lambda_\al>0$.
Among the probability densities on $\setR^d$,
there exists exactly one minimizer $\Blatt_\al\in\pp_2$ of $\rel_\al$,
and it is also the unique minimizer of $\Fishm_\al$.
This follows immediately from Theorem \ref{thm:metric_estimate}.
For $\alpha\in(1/2,1]$, the minimizer $\Blatt_\al=\blatt_\al\lbgd$ is given by
a Barenblatt profile \eqref{eq.barenblatt} of mass $\sfm=1$,
which naturally coincides with the unique stationary solution
of the slow diffusion equation \eqref{eq.slowdiff}.
With the coefficients defined as above, this profile takes the form
\begin{align}
 \label{eq.barenblatt2}
 \blatt_\al(x) = \big( \sfa - \sfb |x|^2 \big)_+^{1/(\alpha-1/2)},
 \quad \sfb = \frac{\alpha-1/2}{\sqrt{2\alpha}} \Lambda,\
 \text{so that}\
 \df\big(\frac {\sqrt{2\alpha}}{2\alpha-1} \blatt_\al^{\alpha-1/2}+
 \frac\Lambda2\,|x|^2\big)=0 .
\end{align}
Here $\sfa>0$ is a parameter chosen to adjust the mass of $\blatt_\al$ to one.
For $\alpha=1/2$, the minimizer $\Blatt_{\frac12,\lambda}$ is a Gaussian measure.

To simplify calculations, we shall use normalized versions of entropy and information,
which are non-negative and vanish exactly in the equilibrium state of \eqref{eq.main}.
Introduce accordingly
\begin{align*}
 \hat{\rel}_\al[\mu] &= \rel_{\alpha,\lambda}[\mu] - \rel_{\alpha,\lambda}[\Blatt_\al], \qquad
 \hat{\Fishm}_\al[\mu] = \Fishm_\al[\mu] - \Fishm_\al[\Blatt_\al],
 \qquad \Blatt_\al=\blatt_\al\lbgd.
\end{align*}
Since $\rel_\al$ is geodesically $\Lambda$-convex in $\pp_2(\Rd)$,
it satisfies the bounds
\begin{align}
 \label{eq.sobolevembed}
 \frac\Lambda2 W_2^2(\mu,\beta_\al)\leq
 \hat{\rel}_\al[\mu] \leq
 \frac 1{2\Lambda}\big|\partial\rel_\al\big|^2(\mu).
\end{align}
As Otto and Villani note, these bounds generalize the
logarithmic Sobolev and Talagrand transportation inequalities \cite{OttoVillani00}.


The equilibration estimate \eqref{eq.l1decay}
is a consequence of the following exponential decay of the entropy.
\begin{theorem}
 \label{prp.prequilibrate}
 Under the assumption of a positive confinement strength $\lambda>0$,
 any solution $\mu(t)$ to problem \eqref{eq.main}\&\eqref{eq.ic}
 constructed in Section \ref{sct.existence} satisfies
 \begin{align}
   \label{eq.entropydecay}
   \frac \Lambda2\,W^2_2(\mu(t),\beta_\al)\leq
   \hat{\rel}_\al[\mu(t)]
   &\leq \exp( - 2\lambda t ) \hat{\rel}_\al[\mu_0], \\
   \label{eq:76}
   \hat{\Fishm}_\al[\mu(t)]&\leq
   \exp( - 2\lambda t ) \hat{\Fishm}_\al[\mu_0].
 \end{align}
\end{theorem}
\begin{proof}
 Since $\rel_\al$ is geodesically $\Lambda$-convex in
 $\pp_2(\Rd)$, it generates a $\Lambda$-flow (see Theorem
 \ref{thm.metricheat}).  Because $\Fishm_\al$ satisfies \eqref{eq:38bis},
 we can apply Theorem \ref{thm:metric_estimate} in the
 complete metric space $X:=\pp_2(\Rd)$ with
 \begin{equation}
   \label{eq:77}
   \V:=\rel_\al,\quad
   \U:=\Fishm_\al,\quad
   \kappa:=\Lambda,\quad
   \theta:=\frac12 \Lambda\,\delta\geq \kappa,\quad
   4\kappa\theta =2\lambda.
 \end{equation}
 By \eqref{eq:20} the discrete solutions $M^n_\sttau$ introduced
 in section \ref{sct.semidiscrete} satisfy
 \begin{equation}
   \label{eq:78}
   (1+2\lambda\tau_n)\hat\rel_\al[M^n_\sttau]\leq
   \hat\rel[M^{n-1}_\sttau],\quad
   (1+2\lambda\tau_n)\hat\Fishm_\al[M^n_\sttau]\leq \hat\Fishm_\al[M^{n-1}_\sttau].
 \end{equation}
 Recalling Remark \ref{rem:exp}, the lower semicontinuity of
 $\rel_\al,\Fishm_\al$, and the pointwise convergence of
 $M_\sttau$ to $\mu$ in $\pp_2(\Rd)$ we conclude both \eqref{eq.entropydecay} and \eqref{eq:76}.
\end{proof}
To conclude estimate \eqref{eq.l1decay} of Theorem \ref{thm.asymptotics},
we invoke a Csiszar-Kullback inequality.
The latter allows to estimate the $L^1$-distance
of a given probability density $\w$ to the Barenblatt profile $\blatt_\al$
in terms of the normalized entropy $\hat{\rel}_\al$.
The following result can be found in \cite[Theorem 30]{CJMTU}.
\begin{lemma}
 \label{lem.csiszar}
 There exists a positive constant $\sfc_\al$ only depending on
 $\alpha$ and $\lambda$ such that for
 every measure $\mu=u\lbgd\in \pp_2^r(\Rd)$
 \begin{align}
   \label{eq.csiszar}
   \| \w - \blatt_\al \|_{L^1(\setR^d)} \leq \sfc_\al
   \, \hat{\rel}_\al[\mu]^{1/2}.
 \end{align}
\end{lemma}
Clearly, \eqref{eq.l1decay} is obtained in combination
of \eqref{eq.csiszar} with \eqref{eq.entropydecay}.

\subsection{Equilibration in $W^{1,2}(\Rd)$ for the thin-film equation}
\label{sct.w12}
In the special case $\alpha=1$ of the thin film equation,
the exponential convergence of the information functionals has the consequence
that the solution $\w(t)$ converges to its steady state in
$W^{1,2}(\setR^d)$,
also exponentially fast.
\begin{theorem}
 Assume $\alpha=1$ and $\lambda>0$.
 Then there exists a constant $\sfC_{\lambda,d}$
 (only depending on $\lambda$ and $d$) such that
 any weak solution $\w(t)$ constructed in the proof of
 Theorem \ref{thm.existence} satisfies
 \begin{align}
   \label{eq.h1decay}
   \big\| \w(t) - \blatt_{1,\lambda} \big\|_{W^{1,2}(\Rd)}^2
   \leq \sfC_{\lambda,d}\hat\fish_{1,\lambda}(u_0) \exp(-2\lambda t).
 \end{align}
 Above, the stationary state
 \begin{align*}
   \blatt_{1,\lambda} (x) =
   \frac{\lambda}{ 8(d+2)}\big( \rho^2 - |x|^2 \big)_+^2 ,\quad
   \text{$\rho=\rho_\lambda>0$ chosen to adjust the mass to unity,}
 \end{align*}
 is the Smyth-Hill profile.
\end{theorem}
\begin{proof}
 For a given probability density $u$,
 define the linear segment $u_\theta=\theta u+(1-\theta) \blatt_{1,\lambda}$ for $0\leq\theta\leq1$.
 Since $\blatt_{1,\lambda}$ is the minimizer of $\fish_{1,\lambda}$,
 \begin{align}
   \label{eq.fishtheta}
   \fish_{1,\lambda}[u_\theta] - \fish_{1,\lambda}[u_0] \geq 0
 \end{align}
 for all $\theta\in[0,1]$.
 Divide \eqref{eq.fishtheta} by $\theta>0$ and calculate the limit $\theta\downarrow0$,
 finding
 \begin{align*}
   \int_\Rd \big( \df u - \df \blatt_{1,\lambda} \big)\cdot \df \blatt_{1,\lambda} \,\dx
   + \frac\lambda2 \int_\Rd |x|^2(u-\blatt_{1,\lambda})\,\dx \geq 0 .
 \end{align*}
 In view of the definition of $\fish_{1,\lambda}$, it follows easily that
 \begin{align}
   \label{eq.h1fromf}
   \frac 12\int_{\setR^d} \big| \df \w - \df \blatt_{1,\lambda}\big|^2 \,\dx
   \leq \fish_{1,\lambda}(\w) -
   \fish_{1,\lambda}(\blatt_{1,\lambda})=
   \Fishm_{1,\lambda}[\mu]-\Fishm_{1,\lambda}[\Blatt_{1,\lambda}]=
   \hat\Fishm_{1,\lambda}[\mu] .
 \end{align}

 The information decay \eqref{eq.informationdecay} thus implies
 \begin{align}
   \label{eq.h1half1}
   \big\| \df \w(t) - \df \blatt_{1,\lambda} \big\|_{L^2(\Rd)}^2
   \leq 2\hat\fish_{1,\lambda}(\w_0)
   \exp\big( - 2\lambda t ).
 \end{align}
 Interpolation of this estimate with the equilibration in $L^1$ from \eqref{eq.l1decay}
 provides the $L^2$-estimate
 \begin{equation}
   \label{eq.h1half2}
   \big\| \w(t) - \blatt_{1,\lambda} \big\|_{L^2(\Rd)}^2
   \leq \hat\ent_{1,\lambda}(\w_0)\,\exp\big(-2\lambda t\big).
 \end{equation}
 Since $\Lambda>0$ and $\hat\Fishm\geq|\partial\rel|^2\geq2\Lambda\hat\rel$,
 the claim \eqref{eq.h1decay} follows by
 combination of \eqref{eq.h1half1} with \eqref{eq.h1half2}.
\end{proof}

\subsection{Equilibration in the absence of confinement}
\label{sct.intermediate}
We shall now consider the gradient flow \eqref{eq.main}
in the limit of vanishing confinement, $\lambda=0$.
In fact,
one may argue that the most natural gradient flow to study in this context
is not \eqref{eq.main} associated to the functional $\fish_{\alpha,\lambda}$,
but the one for the unperturbed information $\fish_{\alpha,0}$.
The corresponding equation reads
\begin{align}
 \label{eq.wholespace}
 \partial_t \www &= -\,\dv\big( \www \df \big[ \www^{\alpha-1} \Delta \www^{\alpha} \big] \big).
\end{align}
Since the unconfined information functional $\Fishm_{\alpha,0}$ enjoys the scaling property \eqref{eq:115},
solutions $\www$ to its gradient flow can be related to solutions $\w$
of a confined gradient flow as follows.
\begin{lemma}
 \label{le:rescaling}
 A function $u\in L^2_{\mathrm{loc}}((0,T);W^{2,2}(\Rd))$ is a weak solution of \eqref{eq.main}
 for $\lambda=1$
 if and only if the function
 \begin{equation}
   \label{eq:98}
   w(t,\cdot)= \Dil{R(t)} u(\log R(t),\cdot)=
   R(t)^{-d} u(\log R(t),R(t)^{-1}\cdot),\quad
   R(t) = \big( 1+ (\delta_\alpha+2)t \big)^{1/(\delta_\alpha+2)},
 \end{equation}
 is a weak solution of \eqref{eq.wholespace}.
\end{lemma}
The proof follows by straight-forward calculations.

An immediate consequence of Lemma \ref{le:rescaling} is that
equation \eqref{eq.wholespace} possesses self-similar solutions.
One is defined by
\begin{align}
 \label{eq.selfsimilar}
 \blatt_{\alpha,0}(t,\cdot) &:= R(t)^{-d}\, \blatt_{\alpha,1}\big(R(t)^{-1}\,\cdot\big)=
 \Dil{R(t)}\blatt_{\alpha,1}(\cdot),
\end{align}
where $\blatt_{\alpha,1}$ is the Barenblatt profile from \eqref{eq.barenblatt2} with unit confinement strength,
and the rescaling factor is given by
\begin{align}
 \label{eq.r}
 R(t) &= \big( 1+ (\delta_\alpha+2)t \big)^{1/(\delta_\alpha+2)}.
\end{align}%
Another consequence of Lemma \ref{le:rescaling} are
the {\em intermediate asymptotics}\nobreakspace of solutions $\www$.
In the limit $t\to\infty$,
any solution $\www$ to \eqref{eq.wholespace} converges to zero in $L^1_{loc}(\setR^d)$,
and does so by approaching the self-similar solution $\www_*$
from \eqref{eq.selfsimilar}.
\begin{corollary}
 \label{cor.intermediate}
 Let $\www$ be generalized minimizing movement induced by $\Fishm_{\alpha,0}$
 (and a fortiori a weak solution to equation \eqref{eq.wholespace}).
 Then $\www$ approaches the self-similar solution $\blatt_{\alpha,0}$ from \eqref{eq.selfsimilar}
 as follows,
 \begin{align}
   \label{eq.intermediateasymptotics}
   \big\| \www(t,\cdot) - b_{\alpha,0}(t,\cdot) \big\|_{L^1(\setR^d)} & \leq C\cdot R(t)^{-1},
 \end{align}
 where $C>0$ depends only on $H_{\alpha,1}(w_0)$ only.
\end{corollary}
%
%
We remark that estimate \eqref{eq.intermediateasymptotics} is a global estimate,
which provides little useful information locally.
Indeed, on any bounded set $M\subset\setR^d$,
the $L^1(M)$-norm of $\www_*(t)$ obviously converges to zero at rate $R(t)^{-d}$,
which is at least as fast as the approximation of $\www_*(t)$ by $\www(t)$ at rate $R(t)^{-1}$.
In order to obtain meaningful local rates, refined asymptotics would need to be studied, as
has been done for the related slow diffusion equation
with $\alpha>1/2$ and $d=1$ by Angenent \cite{Angenent88}, and
with $\alpha<1/2$ but $d\geq 1$ by Denzler, Koch and McCann \cite{DM0} \cite{DKM}.
Globally, however, \eqref{eq.intermediateasymptotics} is a non-trivial statement,
by its equivalence to \eqref{eq.l1decay}.

The goal is to obtain Corollary \ref{cor.intermediate}
as a consequence of estimate \eqref{eq.l1decay} of Theorem \ref{thm.asymptotics}.
The obstacle is that \eqref{eq.l1decay} has not been proven
for general weak solutions to \eqref{eq.main}, but just for those,
which are obtained by means of the minimizing movement variational scheme.
Hence, it remains to be shown that {\em generalized minimizing movements} for the
rescaled and for the unrescaled equation are in correspondence.
To this end, we prove the following (generalized)
discrete counterpart of Lemma \ref{le:rescaling}.
\begin{theorem}
 \label{thm:discrete_rescaling}
 Assume that a functional $\U:\pp_2(\Rd)\to (-\infty,+\infty]$ satisfies
 \begin{equation}
   \label{eq:99}
   \U(\Dil R \mu)=R^{-\delta}\U(\mu)\quad
   \forall\, \mu\in \pp_2(\Rd),\ R>0;\quad
   \Dil R\mu=(R\,\id)_\#\mu.
 \end{equation}
 Let $\bar M\in \pp_2(\Rd),$ $\tau>0$ and
 $R>S>0$ be given.
 A measure $M\in \dom(\U)$ is a minimizer of the functional
 \begin{equation}
   \label{eq:100}
   \mu\mapsto \Phi(\mu;\bar M,\tau,\lambda):=
   \frac 1{2\tau}W_2^2(\mu,\bar M)+\U(\mu)+\frac\lambda2\QMom\mu
 \end{equation}
 if and only if $\Dil R M$ minimizes
 \begin{equation}
   \label{eq:101}
   \begin{gathered}
     \Phi(\mu;\Dil S\bar M,\tilde\tau,\tilde\lambda)= \frac
     1{2\tilde\tau}W_2^2(\mu,\Dil S\bar
     M)+\U(\mu)+\frac{\tilde\lambda}2\QMom\mu,\quad
     \text{with}\\
     \tilde\tau:=\tau S R^{\delta+1},\quad
     \tilde\lambda:=\frac{S(1+\lambda\tau)-R}{\tilde\tau R}.
   \end{gathered}
 \end{equation}
\end{theorem}
In particular, the result applies to $\U=\Fishm_{\alpha,0}$ with $\delta=\delta_\alpha$.
When comparing discrete solutions to \eqref{eq.main} and \eqref{eq.wholespace}, respectively,
$R$ represents the rescaling factor of the current time step,
and $S$ that of the previous one;
Theorem \ref{thm:discrete_rescaling} is applied
with $\tilde\lambda=0$, so that $R=(1+\lambda\tau)S$.
\begin{proof}
 By \eqref{eq:99} and the rescaling property of the Wasserstein distance,
 \begin{equation}
   \label{eq:102}
   W_2^2(\Dil S\mu,\Dil S\nu)=S^2 W_2^2(\mu,\nu)\quad
   \forall\, \mu,\nu\in \pp_2(\Rd),
 \end{equation}
 it is sufficient to prove the theorem in the case $S=1$, when \eqref{eq:101}
 corresponds to
 \begin{equation}
   \label{eq:119}
   \Phi(\mu;\bar M,\tilde\tau,\tilde\lambda)= \frac 1{2\tilde\tau}W_2^2(\mu,\bar M)+\U(\mu)+\frac{\tilde\lambda}2\QMom\mu,\quad
   \text{with}\quad
   \tilde\tau:=\tau R^{\delta+1},\quad
   \tilde\lambda:=\frac{1+\lambda\tau-R}{\tilde\tau R}.
 \end{equation}
 Let us then introduce the functional $\Psi:\dom(\U)\times (0,+\infty)\to\setR$
 \begin{equation}
   \label{eq:103}
   \Psi(\mu,r):=\frac 12 W_2^2(\Dil r\mu,\bar M)+r\,\tau\U(\mu)+
   \frac 12 r(1+\lambda\tau-r) \QMom \mu.
 \end{equation}
 Observe that
 \begin{equation}
   \label{eq:120}
   \tau^{-1}\Psi(\mu,1)=\Phi(\mu;\bar M,\tau,\lambda),\quad
   \big(\tau R^{\delta+1}\big)^{-1}\Psi(\mu,R)=\Phi(\Dil R\mu;\bar M,\tilde\tau,\tilde\lambda)
 \end{equation}
 with $\tilde\tau,\tilde\lambda$ given by \eqref{eq:119}.
 For fixed $\mu\in\pp_2$, the expression $\Psi(\mu,\cdot)$ is linear in its second argument;
 in fact, if $\gamma\in \pp_2(\Rd\times\Rd)$ is an optimal coupling in $\Gamma(\mu,\bar M)$,
 then $(R\id,\id)_\#\gamma$ is an optimal coupling in $\Gamma(\Dil R\mu,\bar M)$.
 Hence,
 \begin{align*}
   W_2^2(\Dil r \mu,\bar M)&=
   \int_{\Rd\times\Rd}|rx-y|^2\,\d\gamma(x,y)=
   \int_{\Rd\times\Rd}\Big(r^2|x|^2-2r x\cdot y+|y|^2\Big)\,\d\gamma(x,y)\\&=
   r^2\QMom\mu-2r\int_{\Rd\times\Rd} x\cdot y\,\d\gamma +\QMom{\bar M},
 \end{align*}
 and
 \begin{displaymath}
   \Psi(\mu,r)=r\Big(-\int_{\Rd\times\Rd} x\cdot y\,\d\gamma +\tau\U(\mu)+
   \frac {1+\lambda\tau}{2} \QMom\mu\Big)+\frac 12\QMom{\bar M}.
 \end{displaymath}
 Differentiating with respect to $r$ we get
 \begin{align*}
   \frac\partial{\partial r}\Psi(\mu,r)&=
   -\int_{\Rd\times\Rd} x\cdot y\,\d\gamma +\tau\U(\mu)+
   \frac {1+\lambda\tau}{2} \QMom\mu\\&=
   \frac 12 W_2^2(\mu,\bar M)+\tau\U(\mu)+\frac{\lambda\tau}{2}\QMom\mu-\frac 12\QMom{\bar M}=
   \Psi(\mu,1)-\frac 12\QMom{\bar M}.
 \end{align*}
 If $M$ is a minimizer of \eqref{eq:100}, then it is also a minimizer of
 $\Psi(\cdot,1)$ and of $\partial_r \Psi(\cdot,r)$.
 It follows that
 \begin{equation}
   \label{eq:105}
   \Psi(M,r)\leq \Psi(\mu,r)\quad\forall\, r>1,\ \mu\in \dom(\U).
 \end{equation}
 In view of the second relation in \eqref{eq:120},
 it follows that $\Dil R M$ is a minimizer of \eqref{eq:101}.
 The same argument shows that if $\tilde M$ is a minimizer in \eqref{eq:101},
 then $\Psi(\Dil {R^{-1}}\tilde M,1)$ coincides with $\Psi(M,1)$,
 and therefore $\Dil {R^{-1}}\tilde M$ is a minimizer in \eqref{eq:100}.
\end{proof}
By means of the previous theorem,
we conclude that generalized minimizing movements indeed enjoy
the same rescaling property stated in Lemma \ref{le:rescaling}.
\begin{corollary}
 \label{cor:rescaling2}
 A continuous curve $\mu\in C^0([0,+\infty);\pp_2(\Rd))$ belongs to
 $GMM(\Fishm_{\alpha,1};\mu_0)$ iff the rescaled curve
 $\nu(t)=\Dil {R(t)}\mu(\log R(t))$ belongs to $GMM(\Fishm_{\alpha,0},\mu_0)$.
\end{corollary}
\begin{proof}
 Let $M_{\sttau}$ be a discrete solution associated to the functional $\Fishm_{\alpha,1}$
 for a given partition $\mathcal P_\sttau$.
 Define a sequence of discrete rescaling $S^n_\sttau$ inductively by
 \begin{equation}
   \label{eq:109}
   S^0_\sttau:=1,\qquad
   S^n_\sttau=(1+\tau_n)S^{n-1}_\sttau\quad
   \text{i.e.\ }\frac {S^n_\sttau-S^{n-1}_\sttau}{\tau_n}=S^{n-1}_\sttau .
 \end{equation}
 We remark that this is the explicit Euler scheme
 for solution of the differential equation
 \begin{equation}
   \label{eq:104}
   \frac {\d}{\d t}S(t)=S(t)
 \end{equation}
 on the partition $\mathcal P_\sttau$.
 Correspondingly,
 we introduce the new family of partitions $\mathcal P_\seeta$,
 induced by
 \begin{equation}
   \label{eq:110}
   \eta_n:=\tau_n S^{n-1}_\sttau \, (S^n_\sttau)^{1+\delta_\alpha},\quad
   s^n_\seeta:=\sum_{m=1}^n \eta_m .
 \end{equation}
 Further,
 we define the piecewise linear time rescaling map $L_\sttau:[0,\infty)\to[0,\infty)$ by
 \begin{equation}
   \label{eq:111}
   L_\sttau(t^n_\sttau):=s^n_\seeta,\quad
   L_\sttau\text{ is linear in }[t^{n-1}_\sttau,t^n_\sttau].
 \end{equation}
 By definition,
 $L_\sttau'(t)=S^{n-1}_\sttau (S^n_\sttau)^{\delta_\alpha+1}$ in $(t^{n-1}_\sttau,t^n_\sttau)$.
 Employing Theorem \ref{thm:discrete_rescaling} with $S:=S^{n-1}_\sttau$ and $R:=S^n_\sttau$,
 it is immediate to check that
 the sequence $\tilde M^n_\seeta$ defined as
 \begin{equation}
   \label{eq:112}
   \tilde M^n_\seeta:= \Dil{S^n}M^n_\sttau
 \end{equation}
 is a discrete solution associated to the unperturbed functional $\Fishm_{\alpha,0}$
 with respect to the rescaled partition $\eeta$.
 The corresponding piecewise constant interpolant $\tilde M_\seeta$ satisfies
 \begin{equation}
   \label{eq:113}
   \tilde M_\seeta(L_\sttau(t))=\Dil{S_\sttau(t)}M_\sttau(t),
 \end{equation}
 with $S_\sttau$ being the piecewise constant interpolant of the $S^n_\sttau$
 with respect to the partition $\mathcal P_\sttau$.

 By definition, $\mu\in GMM(\Fishm_{\al};\mu_0)$ implies that there exists
 a sequence of partitions $\ttau_k$ such that
 $M_{\sttau_k}\to \mu$ locally uniformly in $\pp_2(\Rd)$.
 Moreover, $S_{\sttau_k}(t)\to e^t$ locally uniformly
 and $L_{\sttau_k}(t)\to L(t)=1+(\delta_\alpha+2)^{-1}e^{(\delta_\alpha+2)t}$ by \eqref{eq:110}.
 Thus, up to the extraction of a further subsequence,
 we find that $\tilde M_{\seeta_k}$ converges locally uniformly in $\pp_2(\Rd)$
 to some limit curve $\nu$.
 Again, by definition of the generalized minimizing movements,
 $\nu\in GMM(\Fishm_{\alpha,0};\mu_0)$.
 Passage to the limit in \eqref{eq:113} yields
 \begin{equation}
   \label{eq:121}
   \nu(L(t))=\Dil{e^t}\mu(t)\quad
   \forall\, t>0.
 \end{equation}
 Consequently, $\nu(s)=\Dil{R(s)}\mu(\log R(s))$,
 since $L(t)=s$ means that $t=\log R(s)$.
 The reverse implication can be proved in a similar way.
\end{proof}
Corollary \ref{cor:rescaling2} provides the sought correspondence
between generalized minimizing movements
for equation \eqref{eq.main} and for equation \eqref{eq.wholespace}.
Corollary \ref{cor.intermediate} thus becomes a direct consequence
of estimate \eqref{eq.l1decay} from Theorem \ref{thm.asymptotics}.

\subsection{An alternative approach to entropy decay}
\label{sct.ct}
The following constitutes
an alternative approach to the equilibrium estimate \eqref{eq.entropydecay}.
This alternative approach is purely formal in the range
$\alpha>1/2$, where solutions are not generally smooth,
but might be expected to apply rigorously to the (upper end of)
the range $\alpha \in [(d-2)/2d,1/2]$,  where examples
suggest broader classes of solutions will be smooth (and
possess rapid spatial decay).
We include it here since it provides an enlightening view on the very particular structure
of the gradient flow \eqref{eq.main}.
The key idea is to rewrite the equation \eqref{eq.main} in a way
which is inspired by Carrillo and Toscani \cite{CT2} and
allows us to perform completely explicitly the estimates
that resulted in an abstract way from the theory of gradient flows in the previous section.

Assume that $\w$ is a smooth and strictly positive solution to \eqref{eq.main},
which decays rapidly enough at for $|x|\to\infty$
to justify all subsequent integration by parts.
A tedious but straight-forward calculation reveals that $\w$ satisfies
\begin{align}
 \label{eq.split}
 \partial_t\w &= {\dv}\big( \w \df \FG_\al \big)
 - \frac{2\alpha-1}{2\alpha+1} p \Delta\big(\w^{\alpha+1/2}\Delta \FG_\al \big) -\frac2{2\alpha+1} p \df^2:\big( \w^{\alpha+1/2} \df^2\FG_\al \big)
\end{align}
with the notation
\begin{align*}
 \FG_\al = q \w^{\alpha-1/2} + \frac\lambda2 |x|^2 , \quad
 p = \sqrt{2\alpha} \frac{\Lambda}{\lambda}, \quad
 q = \frac{\sqrt{2\alpha}}{2\alpha-1}
 \frac{\lambda}{\Lambda},
\end{align*}
and $\Lambda$ is as in \eqref{eq.Lambda}.
Notice that the representation for the one-dimensional thin film flow used in \cite{CT2}
is a special case of the general formula \eqref{eq.split} above.

Next, rewrite the definition of $\rel_\al$ in the form
\begin{align*}
 \rel_\al[\mu] = \frac{\Lambda}{\lambda} \int \Big( \frac{q}{\alpha+1/2}\w^{\alpha+1/2} + \frac\lambda2|x|^2\w \Big)\,\dx - C_\alpha .
\end{align*}
It is straight-forward to check that
\begin{align*}
 | \partial\rel_\al |^2 = \Big(\frac\Lambda\lambda\Big)^2 \int \w |\df \FG_\al |^2 \,\dx .
\end{align*}
On the other hand,
the time derivative of the relative entropy along a solution becomes
\begin{align*}
 - \frac{\lambda}{\Lambda}\frac{d}{dt}{\rel}_\al[\mu(t)]
 & =  - \int_{\setR^d} \FG_\al \, \partial_t \w \,\dx \\
 & = \int_{\setR^d} \w \big|\df \FG_\al \big|^2\,\dx \\
 & \qquad + \frac{2\alpha-1}{2\alpha+1} p \int_{\setR^d} \w^{\alpha+1/2} \big| \Delta \FG_\al \big|^2\,\dx
 + \frac2{2\alpha+1} p \int_{\setR^d} \w^{\alpha+1/2} \big\| \df^2 \FG_\al \big\|^2\,\dx.
\end{align*}
Neglect the contribution from the last two terms,
which are non-negative for $\alpha\geq1/2$. (In fact, their sum is non-negative
in the range $\alpha \geq (d-2)/2d$ where $\rel_\al[\mu]$
is $W_2$-geodesically convex.)
Instead, relate the first term of the entropy production to the entropy itself,
observing that
\begin{align*}
 -\frac{d}{dt}\rel_\al[\mu(t)] \geq
 \frac\Lambda\lambda \int_{\setR^d} \w \big| \df \FG_\al \big|^2\,\dx
 = \frac{\lambda}{\Lambda} \, |\partial\rel_\al|^2
 \geq (2\lambda) \, \hat{\rel}_\al,
\end{align*}
employing the convexity estimate\footnote{An
 alternative proof of the second inequality in \eqref{eq.sobolevembed}
 by classical entropy-entropy production methods can be found e.g. in \cite{CJMTU}.}
\eqref{eq.sobolevembed}.
One arrives again at the exponential decay estimate \eqref{eq.entropydecay},
and thus at an alternative proof of Theorem \ref{thm.asymptotics}.
Currently --- lacking estimates on higher regularity of the solution $\w$ ---
it is, however, unclear how to turn this formal argument into a rigorous proof.

\end{document}